\documentclass[times,sort&amp;compress,3p]{elsarticle}
\journal{J. Multivariate Anal.}

\usepackage{lscape,exscale,amsthm, amssymb, booktabs, hyperref}
\usepackage[intlimits]{amsmath}
\usepackage{rawfonts}
\usepackage{latexsym}
\usepackage[cp850]{inputenc}
\usepackage{epsfig}
\usepackage{bbm,enumerate}
\usepackage{soul}
\usepackage{color}
\usepackage{ulem}
\usepackage{stackrel}
\usepackage{multirow}
\usepackage{booktabs}
\usepackage{rotating}
\usepackage{tikz}
\usepackage{float}

\RequirePackage[OT1]{fontenc}
\RequirePackage{amsthm,amsmath}

\newcommand{\bol}[1]{\mbox{\boldmath$#1$}}

\newcommand{\bSigma}{\bol{\Sigma}}

\newcommand{\bm}{\bol{\mu}}

\newcommand{\bx}{\mathbf{X}}

\newcommand{\be}{\mathbf{e}}
\newcommand{\by}{\mathbf{Y}}

\newcommand{\bH}{\mathbf{H}}

\newcommand{\ba}{\mathbf{a}}

\newcommand{\bi}{\mathbf{1}}

\newcommand{\bI}{\mathbf{I}}

\newcommand{\bxi}{\boldsymbol{\xi}}

\newcommand{\bOmega}{\mathbf{\Omega}}

\newcommand{\bS}{\mathbf{S}}
\newcommand{\sx}{\bar{\mathbf{x}}}

\newcommand{\ta}{\alpha}
\newcommand{\tb}{\beta}

\newcommand{\qnullnull}{q_{00}}
\newcommand{\qnn}{q_{nn}}
\newcommand{\qnulln}{q_{0n}}
\newcommand{\snull}{s}
\newcommand{\shatnull}{\hat{s}}
\newcommand{\shatnullc}{\hat{s}_c}
\newcommand{\Rnull}{R}
\newcommand{\Rhatnull}{\hat{R}}

\newcommand{\bZ}{\mathbf{Z}}
\newcommand{\bby}{\bar{\mathbf{y}}}
\newcommand{\bbx}{\bar{\mathbf{x}}}

\newtheorem{theorem}{Theorem}

\newtheorem{lemma}{Lemma}

\newtheorem{remark}{Remark}

\begin{document}
\begin{frontmatter}

\title{Optimal shrinkage estimator for high-dimensional mean vector}

\author[Tmain]{Taras Bodnar\corref{mycorrespondingauthor}}
\cortext[mycorrespondingauthor]{Corresponding author}
\address[Tmain]{Department of Mathematics, Stockholm University, SE-10691 Stockholm, Sweden}
\ead{taras.bodnar@math.su.se}

\author[Omain]{Ostap Okhrin}
\address[Omain]{Chair of Econometrics and Statistics esp. Transportation, Technische Universit\"at Dresden, D-01062 Dresden, Germany}

\author[Nmain]{Nestor Parolya}
\address[Nmain]{Institute of Statistics, Leibniz University of Hannover, D-30167 Hannover, Germany}

\begin{abstract}
 In this paper we derive the optimal linear shrinkage estimator for the high-dimensional mean vector using random matrix theory. The results are obtained under the assumption that both the dimension $p$ and the sample size $n$ tend to infinity in such a way that $p/n \to c\in(0,\infty)$. Under weak conditions imposed on the underlying data generating mechanism, we find the asymptotic equivalents to the optimal shrinkage intensities and estimate them consistently. The proposed nonparametric estimator for the high-dimensional mean vector has a simple structure and is proven to minimize asymptotically, with probability $1$, the quadratic loss when $c\in(0,1)$. When $c\in(1, \infty)$ we modify the estimator by using a feasible estimator for the precision covariance matrix. To this end, an exhaustive simulation study and an application to real data are provided where the proposed estimator is compared with known benchmarks from the literature. It turns out that the existing estimators of the mean vector, including the new proposal, converge to the sample mean vector when the true mean vector has an unbounded Euclidean norm.
\end{abstract}


\end{frontmatter}

\section{Introduction}

Functional data analysis features prominently in modern statistical theory; see  \cite{bongiorno2014contributions, cuevas2014partial, ferraty2006nonparametric, horvath2012inference, ramsay2006functional, ramsay2007applied}. As emphasized in \cite{ferraty2010most, goia2016introduction}, functional data can often be regarded as high-dimensional vectors. This point of view has been explored, e.g., in \cite{aneiros2014variable, aneiros2015partial, aneiros2016sparse, vieu2018}. High-dimensional problems generate a great deal of interest on their own, thanks to the availability of high performance and high storage computers. These advances invite the collection of large amounts of information and stimulate the development of new methods, since classical ones may not be adequate because of specific features of Big Data such as noise accumulation, spurious correlation, heterogeneity and others; see, e.g., \cite{baisil2010}. For a review of challenges in modeling Big Data, see \cite{fanhanlui2014}.

Vast amounts of information are used in various fields, e.g., genomics with hundreds of thousands of microarrays \cite{shalon1996}, neurosciences with high-precision fMRI data of only dozens of persons \cite{issa2013}, financial time series \cite{tsya2010}, etc. In these contexts, one must often deal with the problem of estimating a high-dimensional mean vector when the sample size is much smaller than the dimension. The usual estimator of location under quadratic loss, namely the sample mean vector, has been known for decades to be unsatisfactory because convergence under the quadratic loss function fails when $p$ increases with $n$, even when $p\leq n$. An improved estimator was proposed early on by James and Stein \cite{jamesstein1961} for $p$-dimensional normal random vectors with identity covariance matrix when $p>2$. Baranchik \cite{baranchik1970} then extended this estimator to the case of a diagonal covariance matrix with equal variances. Further research reported in \cite{berger1976, berger1977, fourdrinier2003, gleser1986, lin1973} led to estimators that can be used for arbitrary, unknown covariance matrix for $n\geq p\geq 3$.

More recently, high-dimensional version of the James--Stein type estimator was proposed by Ch\'etelat and Wells \cite{chetelat2012} for $p>n\geq 3$ using an unbiased estimator of the risk difference. They suggested the so-called positive-part type James--Stein estimator and showed through simulations that it dominates the high-dimensional James--Stein estimator under invariant loss. Wang et al. \cite{wang2014} considered an optimal shrinkage estimator towards the unit vector by minimizing expected quadratic loss. However, their estimator is computationally demanding in large dimensions because it involves nontrivial sums. They then suggested that a limiting form of their estimator be used in practice.

This paper contributes to the literature on location parameter estimation by deriving the optimal shrinkage estimator (towards any fixed target) using random matrix theory. The new estimator relies on weaker assumptions than commonly used. We prove its asymptotic normality, derive its limit behavior and show through simulations that it dominates the benchmark methods in terms of both quadratic loss function and computational time.

This paper is structured as follows. The optimal shrinkage estimator for the mean vector is introduced in Section~\ref{sec:2}, where it is shown to be asymptotically equivalent to a nonrandom quantity and Gaussian in the limit. In Section~\ref{sec:3} we provide a bona fide estimator and investigate its asymptotic behavior. Section~\ref{sec:4} discusses benchmark procedures used in the simulation study whose results are reported in Section~\ref{sec:5}. A financial application is presented in Section~\ref{sec:6} and conclusions can be found in Section~\ref{sec:7}. The technical derivations are relegated to the Appendix.

\section{Optimal shrinkage estimator for the mean vector\label{sec:2}}

In this section we construct an optimal shrinkage estimator for the mean vector under large-dimensional asymptotics.
Let $\by_n$ be a $p\times n$ matrix of $p$-dimensional random vectors observed at time points $1,\ldots,n$. The mean vector of each column of $\by_n$ is denoted by $\bm_n$, while $\bSigma_n$ stands for its covariance matrix. Under the large-dimensional asymptotic paradigm, both the dimension $p$ and the sample size $n$ tend to infinity in such a way that $p/n\rightarrow c\in(0,\infty)$ as $n \rightarrow \infty$. For this reason, it is natural to assume, without loss of generality, that $p\equiv p(n)$ is a function of $n$. Later on, we assume that the observation matrix is equal in distribution to
\begin{equation}\label{obs}
 \by_n\overset{d}{=}\bSigma_n^{{1}/{2}}\bx_n+\bm_n\bi_n^\top,
 \end{equation}
where the $p\times n$ matrix $\bx_n$ contains iid random variables with zero mean and unit variance, while $\bi_n$ is the $n$-dimensional vector of $1$s. Only the matrix $\by_n$ is observable. Neither $\bx_{n}$ nor $\bSigma_n$ with $\bm_n$ are known.

It must be noted that the observation matrix $\by_n$ has dependent rows but independent columns. The assumption of independence imposed on its columns can be weakened further to dependent elements of $\bx_n$ by controlling the growth of the number of dependent entries while their joint distribution can be arbitrary; see \cite{friesen2013}. Thus the independence assumption is only used below for technical convenience.

Next, we present the main assumptions which are used in this paper.
\begin{description}
\item[] (A1) There exists $\lambda_0>0$ such that $\lambda_0\le \lambda_{\min}(\bSigma_n)$ uniformly on $p$, where $\lambda_{\min}(\mathbf{A})$ denotes the smallest eigenvalue of the square matrix $\mathbf{A}$. Similarly, the largest eigenvalue of $\mathbf{A}$ is denoted by $\lambda_{\max}(\mathbf{A})$.
\item[] (A2) There exist constants $\gamma >0$, $M_\ell>0$, $M_u>0$ such that $\lim_{p\rightarrow \infty} p^{-\gamma}||\bm_n||^2= M_\ell$ and $\lim_{p\rightarrow \infty} p^{-\gamma}||\bm_0||^2= M_0$ with $0<M_\ell\le M_0,M_\ell \le M_u<\infty$. Here, $\bm_0$ denotes the target mean vector.
\item[] (A3) The elements of the matrix $\bx_n$ have uniformly bounded $2+\varepsilon$ moments for some $\varepsilon>0$.
\end{description}

All of these regularity assumptions are very general and fit many practical situations. Assumption (A1) controls the behavior of the smallest eigenvalue of the population covariance matrix. It is remarkable that no condition is imposed on the largest eigenvalue of $\bSigma_n$, which could also increase to infinity as $p$ becomes larger. This, in particular, allows the application of our results to high-dimensional factor models, which are very popular in economics and finance; see, e.g., \cite{baing2002, baing2008, bodnar_reiss_2016, chamberlain1973, fan2008, fan2013, fan2012}. The assumption on the lower bound of the smallest eigenvalue of $\bSigma_n$ can be avoided via structural assumptions on the true covariance matrix $\bSigma_n$ like tapering; see \cite{cai2010}. In this paper we assume that all eigenvalues are strictly positive. The increase in the norms of the unknown mean vector $\bm_n$ and of the target vector $\bm_0$ are monitored by Assumption (A2), which only requires that they are of the same order. Finally, Assumption (A3) is a technical one and can be relaxed in some cases; see \cite{rubio2012}.

The general linear shrinkage estimator of the mean vector $\bm_n$ is given by
\begin{equation}\label{gse}
\hat{\bm}_{GSE}=\alpha_n\bby_n+\beta_n\bm_0,
\end{equation}
where $\bby_n=n^{-1}\by_n\bi_n$ stands for the sample mean vector. The target vector $\bm_0$ can also be random but independent of the actual information set $\by_n$. The aim is to find the optimal shrinkage intensities which minimize the quadratic loss for a given target vector $\bm_0$ expressed as
$$
L=||\bSigma_n^{-1/2}(\hat{\bm}_{GSE}-\bm_n)||^2=(\hat{\bm}_{GSE}-\bm_n)^\top\bSigma^{-1}_n(\hat{\bm}_{GSE}-\bm_n).
$$
The application of (\ref{gse}) leads to the following optimization problem:
\begin{align*}
&\min_{\ta_n, \tb_n}\ta_n^2\bby_n^\top\bSigma^{-1}_n\bby_n+\tb_n^2\bm_0^\top\bSigma^{-1}_n\bm_0+2\ta_n\tb_n\bby_n^\top\bSigma^{-1}_n\bm_0-2\ta_n\bby_n^\top\bSigma^{-1}_n\bm_n-2\tb_n\bm_n^\top\bSigma^{-1}_n\bm_0.
\end{align*}
Taking the derivatives of $L$ with respect to $\ta_n$ and $\tb_n$ and setting them equal to zero, we get
\begin{align*}
\dfrac{\partial}{\partial\ta_n}\, L &= \ta_n\bby_n^\top\bSigma^{-1}_n\bby_n+\tb_n\bby_n^\top\bSigma^{-1}_n\bm_0-\bby_n^\top\bSigma^{-1}_n\bm_n=0,\\
\dfrac{\partial }{\partial\tb_n}\, L&= \tb_n\bm_0^\top\bSigma^{-1}_n\bm_0+\ta_n\bby_n^\top\bSigma^{-1}_n\bm_0-\bm_n^\top\bSigma^{-1}_n\bm_0=0.
\end{align*}
The Hessian of $L$ is given by
$$
 \bH=\left(
 \begin{array}{cc}
 \bby_n^\top\bSigma^{-1}_n\bby_n&\bby_n^\top\bSigma^{-1}_n\bm_0\\
 \bby_n^\top\bSigma^{-1}_n\bm_0&\bm_0^\top\bSigma^{-1}_n\bm_0
 \end{array}
 \right),
$$
which is a positive definite matrix with probability $1$, since $\bby_n^\top\bSigma^{-1}_n\bby_n>0$ with probability $1$ and
\begin{align}\label{posdef}
 \text{det}(\bH)&=\bby_n^\top\bSigma^{-1}_n\bby_n\bm_0^\top\bSigma^{-1}_n\bm_0-(\bby_n^\top\bSigma^{-1}_n\bm_0)^2 >0
\end{align}
with probability $1$ following the Cauchy--Schwarz inequality applied to the vectors $\bSigma^{-1/2}_n\bby_n$ and $\bSigma^{-1/2}_n\bm_0$.
Thus, the optimal shrinkage intensities are given by
 \begin{align*}
 \alpha_n^*=\alpha_n^*(\bby_n, \bSigma_n, \bm_n, \bm_0) &=  \dfrac{\bby_n^\top\bSigma^{-1}_n\bm_n\bm_0^\top\bSigma^{-1}_n\bm_0-\bm_n^\top\bSigma^{-1}_n\bm_0\bby_n^\top\bSigma^{-1}_n\bm_0}{\bby_n^\top\bSigma^{-1}_n\bby_n\bm_0^\top\bSigma^{-1}_n\bm_0-(\bby_n^\top\bSigma^{-1}_n\bm_0)^2},\\
 \beta_n^*=\beta_n^*(\bby_n, \bSigma_n, \bm_n, \bm_0)&= \dfrac{\bby_n^\top\bSigma^{-1}_n\bby_n\bm_n^\top\bSigma^{-1}_n\bm_0-\bby_n^\top\bSigma^{-1}_n\bm_0\bby_n^\top\bSigma^{-1}_n\bm_n}{\bby_n^\top\bSigma^{-1}_n\bby_n\bm_0^\top\bSigma^{-1}_n\bm_0-(\bby_n^\top\bSigma^{-1}_n\bm_0)^2}.
 \end{align*}

In Theorem~\ref{th1}, we show that the optimal shrinkage intensities $\alpha^*_n$ and $\beta^*_n$ are almost surely asymptotically equivalent to nonrandom quantities $\alpha^*$ and $\beta^*$ under the large-dimensional asymptotics $p/n\rightarrow c\in(0,\infty)$ as $n\rightarrow \infty$.

\begin{theorem}\label{th1}
Assume that (A1)--(A3) hold. Then
$|\alpha^*_n-\alpha^* | \stackrel{\rm a.s.}{\to }0$ and $ |\beta^*_n-\beta^* |\stackrel{\rm a.s.}{\to }0$
for $p/n\rightarrow c>0$ as $n\rightarrow\infty$, where
\begin{align}\label{a_as}
\alpha^*&=  \alpha^*(\bSigma_n, \bm_n, \bm_0)=\frac{\bm_n^\top\bSigma^{-1}_n\bm_n \bm_0^\top\bSigma^{-1}_n\bm_0-  (\bm_n^\top\bSigma^{-1}_n\bm_0)^2}
{(c+ \bm_n^\top\bSigma^{-1}_n\bm_n) \bm_0^\top\bSigma^{-1}_n\bm_0-(\bm_n^\top\bSigma^{-1}_n\bm_0)^2},\\
\label{b_as}
\beta^*&= \beta^*(\bSigma_n, \bm_n, \bm_0)=(1-\alpha^*) \, \dfrac{\bm_n^\top\bSigma^{-1}_n\bm_0}{\bm_0^\top\bSigma^{-1}_n\bm_0} .
\end{align}
\end{theorem}

Note that $\ta^*\in(0,1)$ due to inequality \eqref{posdef}. Furthermore, using the results of Theorem~\ref{th1}, we are able to estimate $\ta^*$ and $\tb^*$ consistently at least for $c\in (0,1)$, which is shown in Theorem~\ref{th3} below.

It is remarkable that the proposed procedure is very different to the one suggested by Wang et al. \cite{wang2014}, who minimized the  expected quadratic loss and estimated optimal shrinkage intensities. While they found the estimators for the optimal shrinkage intensities which converge in probability, our aim is to construct consistent estimators which converge almost surely. It is worth pointing out the following remark.

\begin{remark}
{\rm
The technical assumption {\rm (A2)} with the same $\gamma$s is key. If this condition fails, i.e., there exist $\gamma_1$ and $\gamma_2$ such that $\lim_{p\rightarrow \infty} p^{-\gamma_1}||\bm_n||^2= M_\ell$ and $\lim_{p\rightarrow \infty} p^{-\gamma_2}||\bm_0||^2= M_0$ with $0<M_\ell\le M_0,M_\ell \le M_u<\infty$, then from the expression of $\beta_n^*$ we have that the rate for the numerator is $\gamma_1+(\gamma_1+\gamma_2)/2$ whereas the rate for the denominator is $\gamma_1+\gamma_2$. Consequently, if $\gamma_1 \neq \gamma_2$ then}
\[
\beta_n^* \stackrel{\rm a.s.}{\to }
\left\{
\begin{array}{cc}
  0~~{\rm for}~ n^{-1}p\rightarrow c>0~ {\rm as} ~n\rightarrow\infty & {\rm if}~ \gamma_1 < \gamma_2 ,\\
  \infty~~{\rm for}~ n^{-1}p\rightarrow c>0~ {\rm as} ~n\rightarrow\infty & {\rm if}~ \gamma_1 > \gamma_2.
 \end{array}
\right.
\]
\end{remark}

Let $\tilde{c}= p^{-\gamma}c$, $q_{ij} = p^{-\gamma}\bm_i^\top\bSigma^{-1}_n\bm_j$, for $i, j\in\{0, n\}$ and $d = q_{00}q_{nn}-q_{0n}^2$. In Theorem~\ref{th2}, we prove that $\ta^*_n$ and $\tb^*_n$ are asymptotically Gaussian under the large-dimensional asymptotic regime. It must be noted that the  normality assumption is used only in Theorems \ref{th2} and \ref{th4}, whereas the existence of the second (fourth) moments is only required for the rest of our results.

\begin{theorem}\label{th2}
		Assume (A1)--(A2) and let the elements of $\bx_n$ be standard normally distributed. Then
$$
		\sqrt{p^{\gamma}n} \, \sigma^{-1}_{\alpha} \, (\alpha_n^*-\alpha^*) \rightsquigarrow\mathcal{N}(0,1),\quad
		\sqrt{p^{\gamma}n} \, \sigma^{-1}_{\beta} \, (\beta_n^*-\beta^*)\rightsquigarrow\mathcal{N}(0,1),
$$
		for $p/n\rightarrow c>0$ as $n\rightarrow\infty$, where
$		\sigma^2_{\alpha} =  \{(\tilde{c}q_{00}-d)^2q_{00}d + \tilde{c}d^2q_{00}^2\}(\tilde{c}q_{00}+d)^{-4}$ and
		\begin{multline*}
		\sigma^2_{\beta}=  \frac{1}{(\tilde{c}q_{00}+d)^4}\{(d-\tilde{c}q_{00})^2q_{0n}^2q_{nn}+(\tilde{c}q_{0n}^2-\tilde{c}d-dq_{nn})^2q_{00}+\tilde{c}d^2q_{0n}^2
        +2(\tilde{c}q_{0n}^2-\tilde{c}d-dq_{nn})(d-\tilde{c}q_{00})q_{0n}^2\}.
		\end{multline*}
\end{theorem}

\section{Bona fide estimator\label{sec:3}}

This section presents consistent estimators for $\alpha^*$ and $\beta^*$, i.e., for the deterministic equivalent quantities to the optimal shrinkage estimators $\alpha_n^*$ and $\beta_n^*$, which we denote by $\hat{\alpha}^*$ and $\hat{\beta}^*$. This procedure allows us to construct bona fide estimators for the unknown shrinkage intensities. Using recent results from random matrix theory, we further prove that $\hat{\alpha}^*$ and $\hat{\beta}^*$ are consistent and asymptotically normally distributed.

Let
$$
 \bS_n=n^{-1} (\by_n-\bby_n\bi_n^\top ) (\by_n-\bby_n\bi_n^\top )^{\top}=n^{-1}\by_n\by_n^\top-\bby_n\bby_n^\top
$$
be the sample covariance matrix. In Theorem~\ref{th3} below we present consistent estimators for $\alpha^*$ and $\beta^*$ under  large-dimensional asymptotics.

\begin{theorem}\label{th3}
Assume (A1)--(A2) and let the elements of $\bx_n$ possess uniformly bounded $4+\varepsilon$ moments with $\varepsilon>0$. Then the consistent estimators for $\alpha^*$ and $\beta^*$ are given by
\begin{align}\label{a_con}
\hat{\alpha}^*&=\hat{\alpha}^*(\bby_n, \bS_n, \bm_0)=
\frac{ \{\bby_n^\top \bS^{-1}_n\bby_n - {p}/{(n-p)} \}\bm_0^\top \bS^{-1}_n\bm_0-(\bby_n^\top \bS^{-1}_n\bm_0)^2}{\bby_n^\top \bS^{-1}_n\bby_n \bm_0^\top \bS^{-1}_n\bm_0-(\bby_n^\top \bS^{-1}_n\bm_0)^2},\\
\label{b_con}
\hat{\beta}^*&=\hat{\beta}^*(\bby_n, \bS_n, \bm_0)=(1-\hat{\alpha}^*) \, \dfrac{\bby_n^\top\bS^{-1}_n\bm_0}{\bm_0^\top\bS^{-1}_n\bm_0}.
\end{align}
\end{theorem}

Next, we prove that the consistent estimators for the shrinkage intensities are asymptotically normally distributed. This result is investigated under an additional condition imposed on the distribution of the entries of $\bx_n$, which are assumed to be standard normally distributed.

\begin{theorem}\label{th4}
Assume (A1)--(A2) and let the elements of $\bx_n$ be standard normally distributed. Then
\begin{align*}
\sqrt{n}\, \bOmega^{-1/2}\left(\begin{array}{c}
  \hat{\alpha}^*-\alpha^* \\
  \hat{\beta}^*-\beta^*\\
          \end{array}
        \right)\rightsquigarrow  \mathcal{N}\left[\left(
                                                   \begin{array}{c}
                                                     0 \\
                                                     0 \\
                                                   \end{array}
                                                 \right),\mathbf{I} \right]
\end{align*}
where
\begin{align*}
\bOmega&= \left(\begin{array}{cc}  {c^2\sigma_s^2}/{(c+\snull)^4} & {c^2 \sigma_s^2\Rnull }/{(c+\snull)^4}   \\
{c^2 \sigma_s^2\Rnull}/{(c+\snull)^4}  & {c^2 \sigma_s^2\Rnull ^2}/{(c+\snull)^4} +  {c^2}{(c+\snull)^{-2}} \{1+{(\snull+c)}/{(1-c)}\}/{\bm_0^\top \bSigma^{-1}_n\bm_0} \\\end{array}\right)
\end{align*}
and
\[\snull=\bm_n^\top \bSigma^{-1}_n\bm_n-\frac{(\bm_0^\top \bSigma^{-1}_n\bm_n)^2}{\bm_0^\top \bSigma^{-1}_n\bm_0},\quad \Rnull =\frac{\bm_0^\top \bSigma^{-1}_n\bm_n}{\bm_0^\top \bSigma^{-1}_n\bm_0} , \quad
	\sigma_s^2= 2\left(c+2\snull\right)+\frac{2}{1-c}\left(c+\snull\right)^2.
\]
\end{theorem}

The results of Theorem~\ref{th4} can be used to construct an asymptotic test on both $\alpha^*$ and $\beta^*$. For example, testing the null hypothesis $\mathcal{H}_0:~\alpha^*=1$ will imply that in terms of the quadratic loss, the best estimator for the mean vector  is the sample mean vector because $\beta^*=0$ as soon as $\alpha^*=1$.

The bona fide optimal shrinkage estimator for the mean vector in the case of $c<1$ is constructed by
\begin{equation}\label{olse1}
\hat{\bm}_{OLSE}=\hat{\alpha}^*\bby_n+\hat{\beta}^*\bm_0,
\end{equation}
where the optimal shrinkage intensities are given by \eqref{a_con} and \eqref{b_con}, respectively. This estimator has almost surely smallest quadratic loss under the large-dimensional asymptotics. We refer to it as the Optimal Linear Shrinkage Estimator (OLSE) for the high-dimensional mean vector. It is obvious that the OLSE estimator \eqref{olse1} dominates the sample estimator in terms of minimum quadratic loss uniformly if both $p$ and $n$ tend to infinity and $p/n\rightarrow c<1$.

For $c>1$ the sample covariance matrix is no longer invertible and we need other techniques to estimate the unknown quantities given in \eqref{a_as} and \eqref{b_as}. Here, we apply the generalized inverse of the sample covariance matrix $\bS_n$. Particularly, we use the following generalized inverse of the sample covariance matrix $\bS_n$:
$$
\bS_n^-=\bSigma_n^{-1/2} ( \bx_n\bx_n^\top/n - \sx_n\sx_n^\top )^+\bSigma_n^{-1/2},
$$
where $^+$ denotes the Moore--Penrose inverse. It can be shown that $\bS_n^-$ is a generalized inverse of $\bS_n$ satisfying $\bS_n^-\bS_n\bS_n^-=\bS_n^-$ and $\bS_n\bS_n^-\bS_n=\bS_n$. However, $\bS_n^-$ is not exactly equal to the Moore--Penrose inverse because it does not satisfy the conditions $(\bS_n^-\bS_n)^\top=\bS_n^-\bS_n$ and $(\bS_n\bS_n^-)^\top=\bS_n\bS_n^-$. When $c<1$, the generalized inverse $\bS_n^-$ coincides with the usual inverse $\bS_n^{-1}$. Moreover, if $\bSigma_n$ is a multiple of identity matrix then $\bS_n^-$ is equal to the Moore--Penrose inverse $\bS_n^+$. So, it could be expected that if $\bSigma_n$ is a sparse matrix, then both the inverses are very close. This conjecture is not shown here; it is left for future research.

In Theorem~\ref{th5} below we present the consistent estimators for $\alpha^*$ and $\beta^*$ under large-dimensional asymptotics in the case of $c>1$ utilizing the generalized inverse $\bS_n^-$.

\begin{theorem}\label{th5}
Assume (A1)--(A2) and let the elements of $\bx_n$ possess uniformly bounded $4+\varepsilon$ moments with $\varepsilon>0$. Let $p/n\rightarrow c\in(1, \infty)$ for $n\rightarrow\infty$. Then  consistent estimators for $\alpha^*$ and $\beta^*$ are given by
\begin{align*}
\hat{\alpha}^*&= \hat{\alpha}^*(\bby_n, \bS_n, \bm_0)=
\frac{\{\bby_n^\top \bS^{-}_n\bby_n -(p/n-1)^{-1} \} \bm_0^\top \bS^{-}_n\bm_0-(\bby_n^\top \bS^{-}_n\bm_0)^2}{\bby_n^\top \bS^{-}_n\bby_n \bm_0^\top \bS^{-}_n\bm_0-(\bby_n^\top \bS^{-}_n\bm_0)^2},\\
\hat{\beta}^*&= \hat{\beta}^*(\bby_n, \bS_n, \bm_0)=(1-\hat{\alpha}^*) \, \dfrac{\bby_n^\top\bS^{-}_n\bm_0}{\bm_0^\top\bS^{-}_n\bm_0}.
\end{align*}
\end{theorem}

Because $\bS_n^-$ depends on unknown quantities, we will approximate it by the Moore--Penrose inverse $\bS_n^+$. The asymptotic properties of the Moore--Penrose inverse under high-dimensional settings are investigated in \cite{bodnar_dette_parolya_2016}. It is worth mentioning that changing $\bS^-$ to $\bS^+$ does not in general lead to the optimal shrinkage estimator. This is in contrast to the case $c<1$, where it obviously holds that $\bS^{-1}=\bS^-=\bS^+$. Consequently, the approach suggested for $p >n$ is only suboptimal, but nevertheless it dominates in most cases the existent estimators for the high-dimensional mean vector given in the literature. This fact is justified in the next sections via a simulation study and an empirical illustration.

\subsection{Choice of $\bm_0$}
\label{Sec:ChoiceMu}
An important issue is the choice of the nonrandom target vector $\bm_0$ which should satisfy Assumption (A2). This depends on the underlying data because the choice of the target vector is equivalent to the choice of the hyperparameter for the prior distribution of $\bm_n$. This problem is well-known in Bayesian statistics. Different priors lead to different results. So it is crucial to choose the one which works satisfactory for most cases. The naive choice is $\bm_0=p^{(\gamma-1)/2} \bi$, where $\bi$ is the $p$-dimensional vector of $1$s. Obviously, the perfect $\bm_0$ is the true mean vector $\bm_n$; we then get from Theorem~\ref{th1} that $\alpha^*=0$ and $\beta^*=1$. Hence, the proposed optimal shrinkage estimator will tend to the true mean almost surely as both $p$ and $n$ tend to infinity; setting $\bm_0$ close to $\bm_n$ will improve the resulting shrinkage estimator.

Another quantity present in the naive target $\bm_0=p^{(\gamma-1)/2} \bi$ is $\gamma$, which measures how quickly the Euclidean norm of the true mean vector goes to infinity. As can be seen from the simulation study, if $\gamma \neq 0$, then all of the considered estimators including our proposal converge to the sample mean vector.

\section{Benchmark methods\label{sec:4}}

This section introduces approaches used as benchmarks in the simulation study of Section~\ref{sec:5}. The most commonly used estimator for the mean vector in the literature is the sample mean vector expressed as
$	\bby_n= \by_n\bi_n/n$.
Although this estimator is known to be not converging under quadratic loss for high-dimensional data, we nevertheless use the sample mean in our comparison study.

The original James--Stein estimator was expressed in the form
\begin{equation*}
\hat\bm_{n,JSN}=\left(1-\frac{p-2}{n\bby_n^\top\bby_n}\right) \bby_n
\end{equation*}
for $\bSigma_n=\mathbf{I}_p$ and $n>p>2$.
In the comparison study we make use of a modified version of this estimator given by
\begin{equation*}
	\widehat{\bm}_{n,JS}=\left\{1-\frac{(p-2)/(n-p-3)}{\bby_n^\top \widetilde\bS_n^{-1} \bby_n}\right\} \bby_n
\end{equation*}
for $c<1$ and an estimator $\widetilde\bS_n\sim \mathcal{W}_p(n,\bSigma_n)$ of the covariance matrix $\bSigma_n$ with $n\geq p\geq 3$. When $p>n\geq 3$, we compare our estimator with those proposed by Ch\'etelat and Wells \cite{chetelat2012}, who defined a Baranchik type estimator as
\begin{equation*}
	\widehat{\bm}_{n,JS(p>n)}=\left( \bI_p-\frac{a\widetilde\bS_n\widetilde\bS_n^+}{\bby_n^\top \widetilde\bS_n^+ \bby_n}\right) \bby_n
\end{equation*}
with $0\leq a\leq  {2(n-2)}/{(p-n+3)}$ and $\widetilde\bS_n^+$ the Moore--Penrose inverse of $\widetilde\bS_n$. In our study we set $a =  {2(n-2)}/{(p-n+3)}$. As shown by the simulation study in \cite{chetelat2012}, the so-called positive-part type James--Stein estimator of the form

\begin{equation*}
        \widehat{\bm}_{n,JS+}=(\bI_p + \widetilde\bS_n\widetilde\bS_n^+)\bby_n+\left\{1-\frac{(n-2)/(p-n+3)}{\bby_n^\top \widetilde\bS_n^+ \bby_n}\right\}_+ \widetilde\bS_n\widetilde\bS_n^+\bby_n,
    \end{equation*}
with $b_+ = \max(b, 0)$ dominates $\widehat{\bm}_{n,JS(p>n)}$ under the invariant loss.

Another benchmark estimator is taken from \cite{wang2014}. It is a shrinkage estimator with unit target vector and shrinkage coefficients found by the minimization of the expected quadratic loss. This shrinkage estimator is given by
\begin{equation}\label{bm_nW}
	\widehat{\bm}_{n, W} = \frac{Z_{1,n} - Z_{4,n} }{Z_{1,n} + Z_{2,n}  Z_{4,n}}  \bby_n + \frac{Z_{2,n}}{Z_{1,n} + Z_{2,n}  Z_{4,n}} Z_{3,n} \bi_n,
\end{equation}
with
\begin{align*}
Z_{1,n} & = \frac{1}{p(n-1)} \sum\limits_{i \ne j} \by_{n,i}^\top \widetilde\bS_n^+ \by_{n,j},
\quad Z_{2,n}  = \frac{1}{np} \left(\sum\limits_{k=1}^{n} \by_{n,k}^\top \widetilde\bS_n^+ \by_{n,k}  - \frac{1}{n-1}\sum\limits_{i \ne j} \by_{n,i}^\top \widetilde\bS_n^+ \by_{n,j} \right) ,\\
Z_{3,n} & = \frac{1}{n \bi_n^\top \widetilde\bS_n^+ \bi_n} \sum\limits_{k=1}^n \bi_n^\top \widetilde\bS_n^+ \by_{n,k},
\quad Z_{4,n}  = \frac{1}{p(n-1)\bi_n^\top \widetilde\bS_n^+ \bi_n} \sum\limits_{i \ne j} \bi_n^\top \widetilde\bS_n^+ \by_{n,i} \by_{n,j}^\top \widetilde\bS_n^+ \bi_n,
\end{align*}
for $ p/n>1$. The estimator \eqref{bm_nW} has a computationally complicated form because of the double sum over $p$ and~$n$, being therefore very time consuming for large dimensions and large sample sizes. That is why, in practice, its asymptotic counterpart is considered; see, e.g., \cite{wang2014}.

\section{Finite-sample performance\label{sec:5}}

This section provides an extensive simulation study, to test the validity of Theorems~\ref{th1}, \ref{th2}, \ref{th4}, and \ref{th5} as well as to compare the quality of the proposed OLSE estimator with the considered benchmark methods. As discussed in Section~\ref{Sec:ChoiceMu}, $\gamma$ controls the speed of divergence of the Euclidean norms of $\bm_n$ and $\bm_0$ to infinity. Therefore, we consider two extreme scenarios, namely $\gamma = 0$ and $\gamma = 1$.

When $\gamma = 0$, the true mean vector $\bm_n$ and the target mean vector $\bm_0$ have been independently simulated from the uniform distribution on $[-p^{-1/2}, p^{-1/2}]$. This scenario keeps the norms of both $\bm_n$ and $\bm_0$ bounded. In contrast, only very little prior information about $\bm_n$ is used in the selection of $\bm_0$. This choice of $\bm_n$ is also motivated by the empirical illustration of Section~\ref{sec:6}, where the obtained theoretical results are applied to financial data consisting of asset returns which usually possess small expected values.

When $\gamma = 1$, the true mean vector $\bm_n$ contains randomly interchanging values $1$ and $-1$, while the $\bm_0$ is set to the unit vector. This scenario forces the norms of both mean vectors to explode with the dimension $p$ and it uses no information about the true mean vector $\bm_n$ when $\bm_0$ is chosen.

For all of the considered setups we also construct a bona fide shrinkage estimator with $\bm_0=\bm_n$ which corresponds to the case when the optimal shrinkage estimator converges to $\bm_n$ almost surely; see Section~\ref{Sec:ChoiceMu}. The corresponding results are depicted in all figures as the solid red line. The eigenvalues of the true covariance matrix $\bSigma_n$ are chosen in the following way: 20\% of all eigenvalues are equal to $1$, 40\% of them equal $3$ and the rest 40\% are 10. The eigenvectors are extracted via the QR decomposition of a standard normally distributed random matrix. The obtained eigenvalues and eigenvectors are finally wrapped via the singular value decomposition into a covariance matrix $\bSigma_n$. We further note that the quantities $\bm_n$, $\bm_0$, and $\bSigma_n$ are fixed in each simulation study, i.e., the same values of $\bm_n$, $\bm_0$, and $\bSigma_n$ are used to generate samples which follow model \eqref{obs}.

First, we considered the finite-sample behavior of the optimal shrinkage coefficients $\alpha_n^*$ and $\beta_n^*$ and their bone fide estimators $\hat\alpha^*$ and $\hat\beta^*$. We compared them to the corresponding asymptotic distributions presented in Theorems \ref{th2} and~\ref{th4}. For $p \in \{20, 100, 250, 500\}$ and $c \in \{0.5, 0.9, 2.0\}$, we simulated $\bm_n$, $\bm_0$, $\bSigma_n$ as described above, considered only the case $\gamma = 0$, and drew the columns of $\mathbf{Y}_n$ from $\mathcal{N}(\bm_n, \bSigma_n)$ in each simulation run. Then, $\alpha_n^*$ and $\beta_n^*$ were calculated and standardized by the corresponding asymptotic variances from Theorem~\ref{th2} for $c\in\{0.5, 0.9, 2.0\}$. For $c\in\{0.5, 0.9\}$, we also computed $\hat\alpha^*$ and $\hat\beta^*$ and standardized them by their asymptotic variances as presented in Theorem~\ref{th4}. The procedure was repeated $N = 1000$ times resulting in $N$ values of $\alpha_n^*$, $\beta_n^*$, $\hat\alpha^*$, and $\hat\beta^*$ for chosen $p$ and~$c$.

\begin{figure}[!ht]
		\begin{center}
		\includegraphics[scale=0.5]{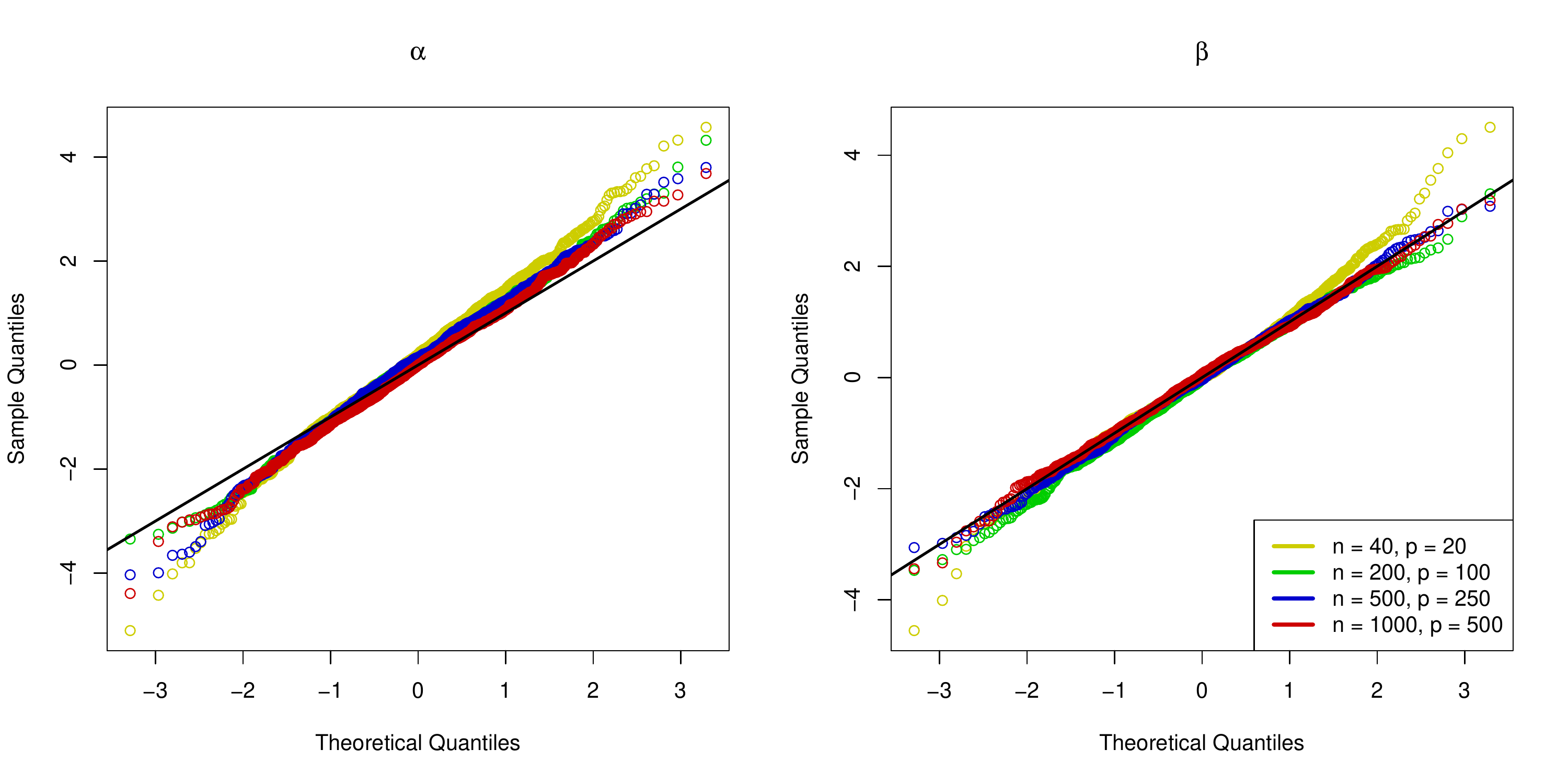}
		\includegraphics[scale=0.5]{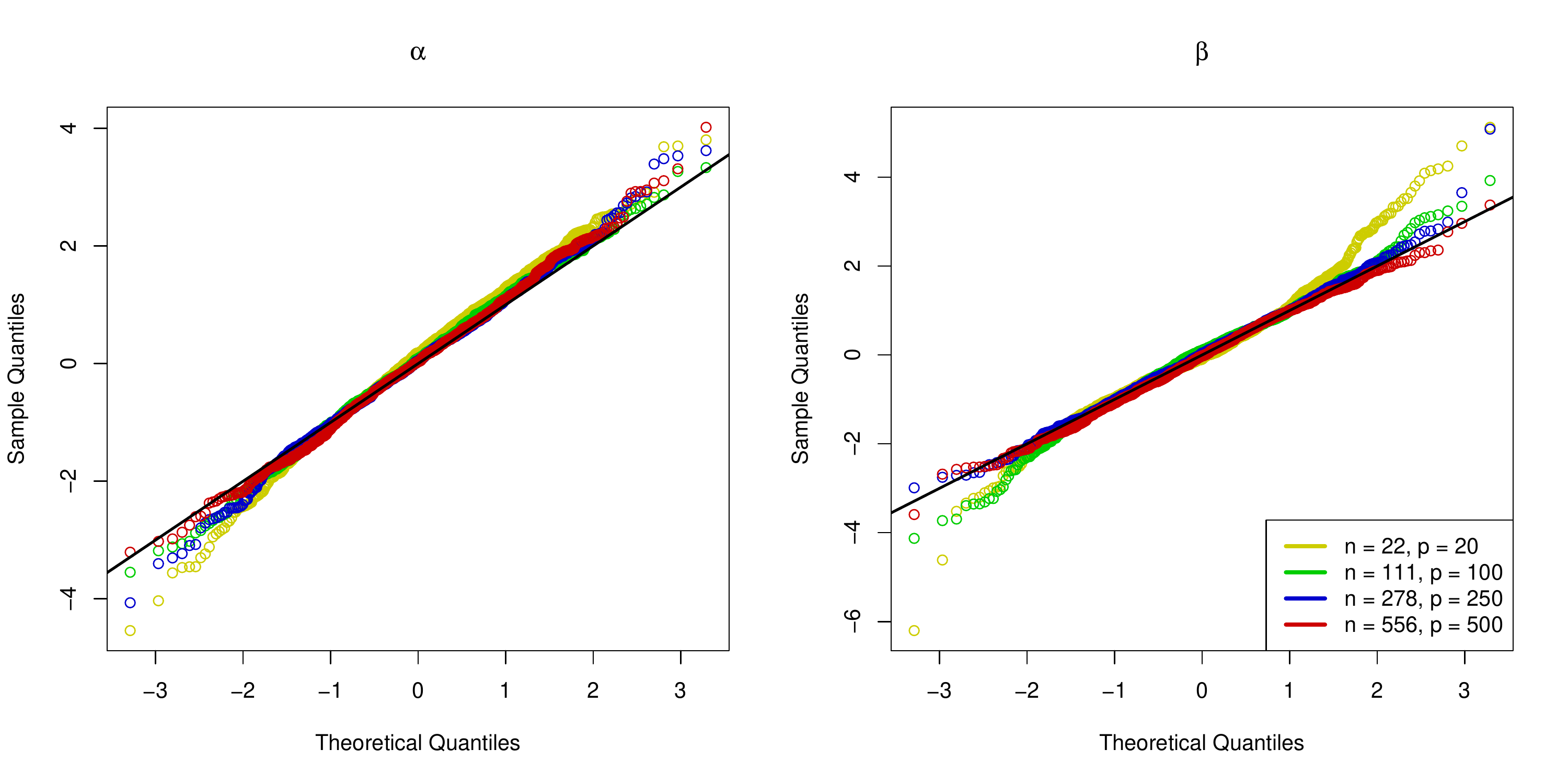}
		\includegraphics[scale=0.5]{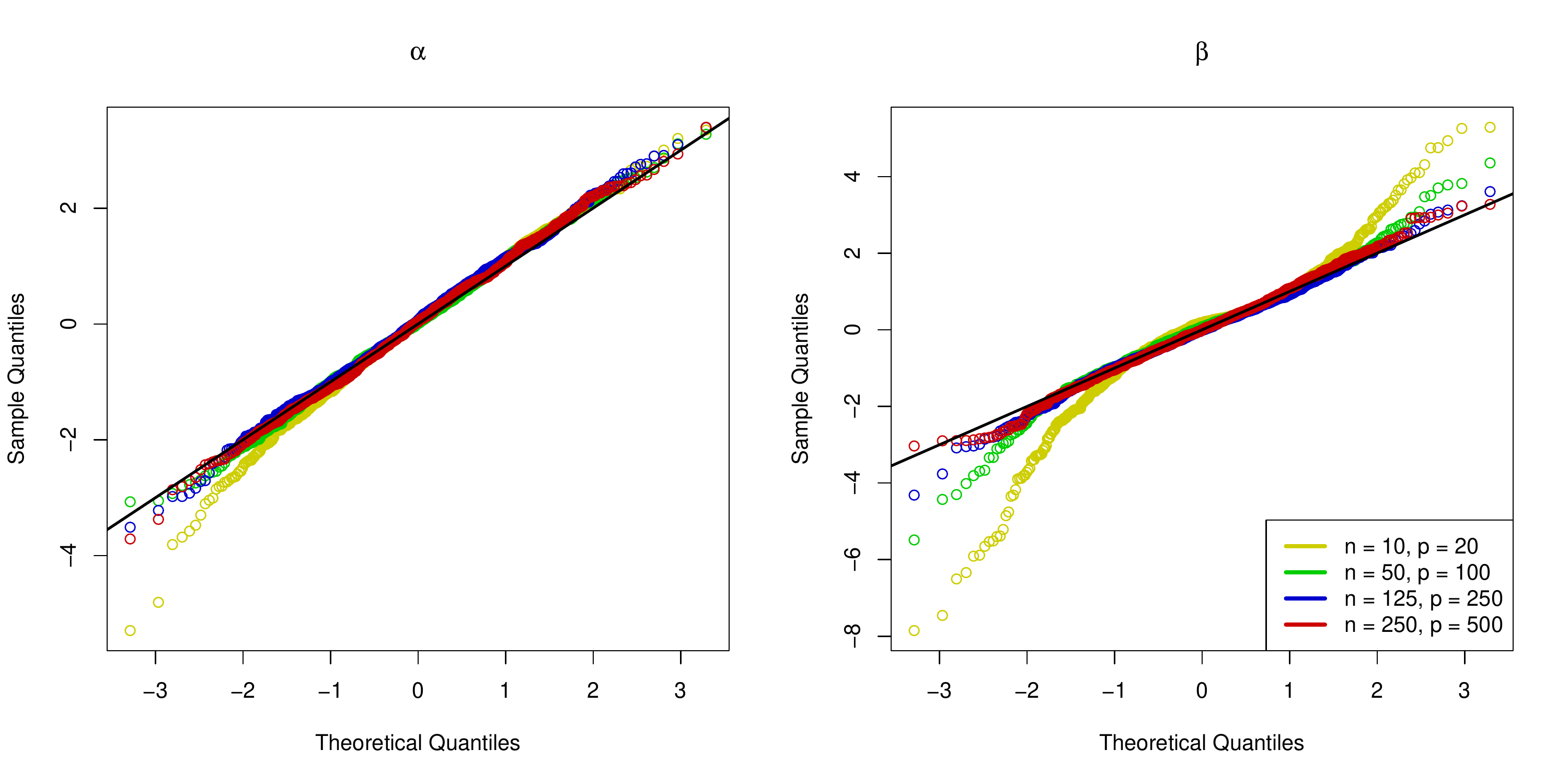}
		\caption{Q-Q plots for $\alpha^*_n$ and $\beta^*_n$ in the case of $\gamma = 0$. We set $c = 0.5$ (top panel), $c = 0.9$ (middle panel) and $c = 2.0$ (bottom panel). }\label{fig:norm:or}
		\end{center}
	\end{figure}

\begin{figure}[!ht]
		\begin{center}
		\includegraphics[scale=0.5]{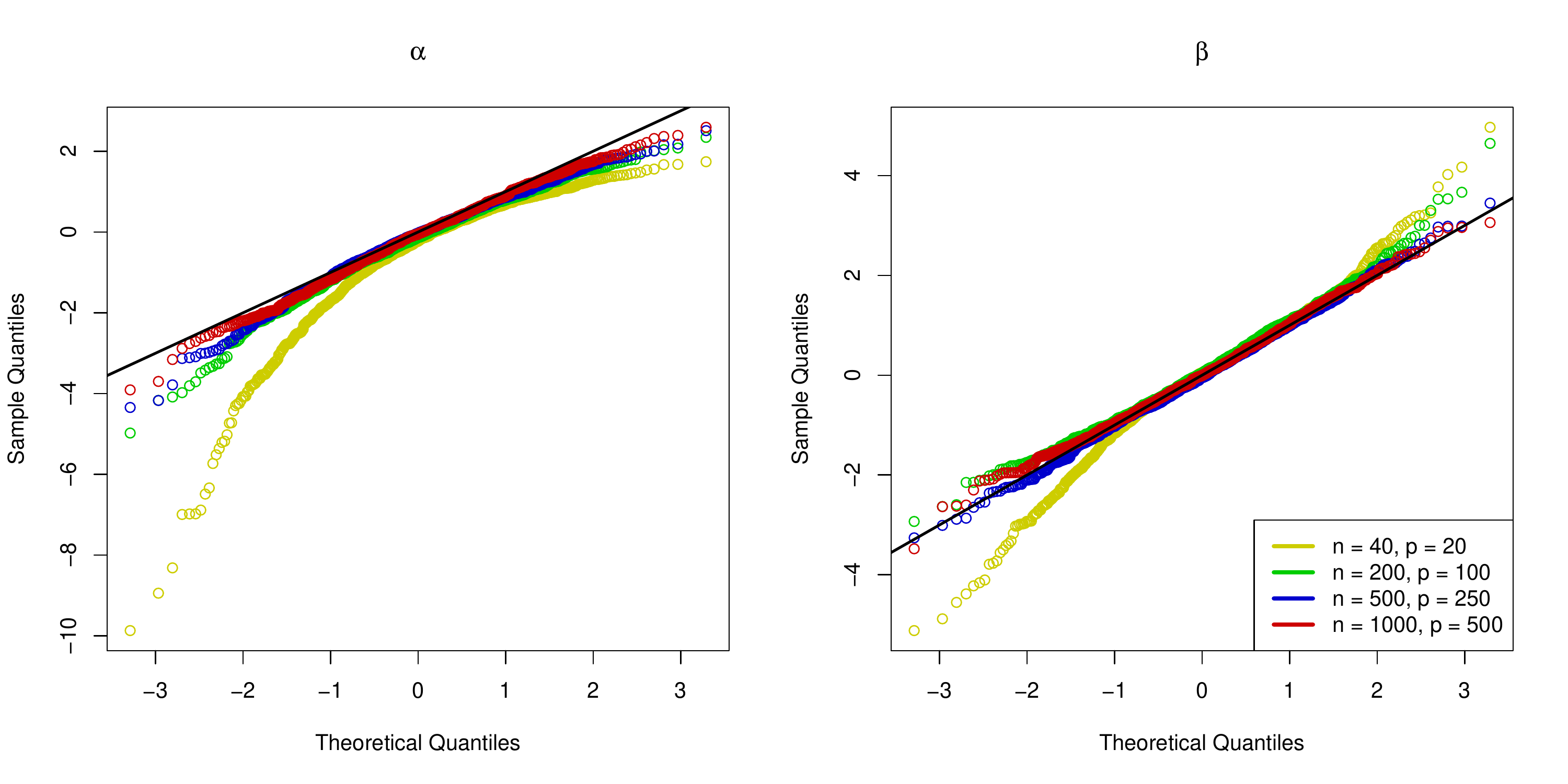}
		\includegraphics[scale=0.5]{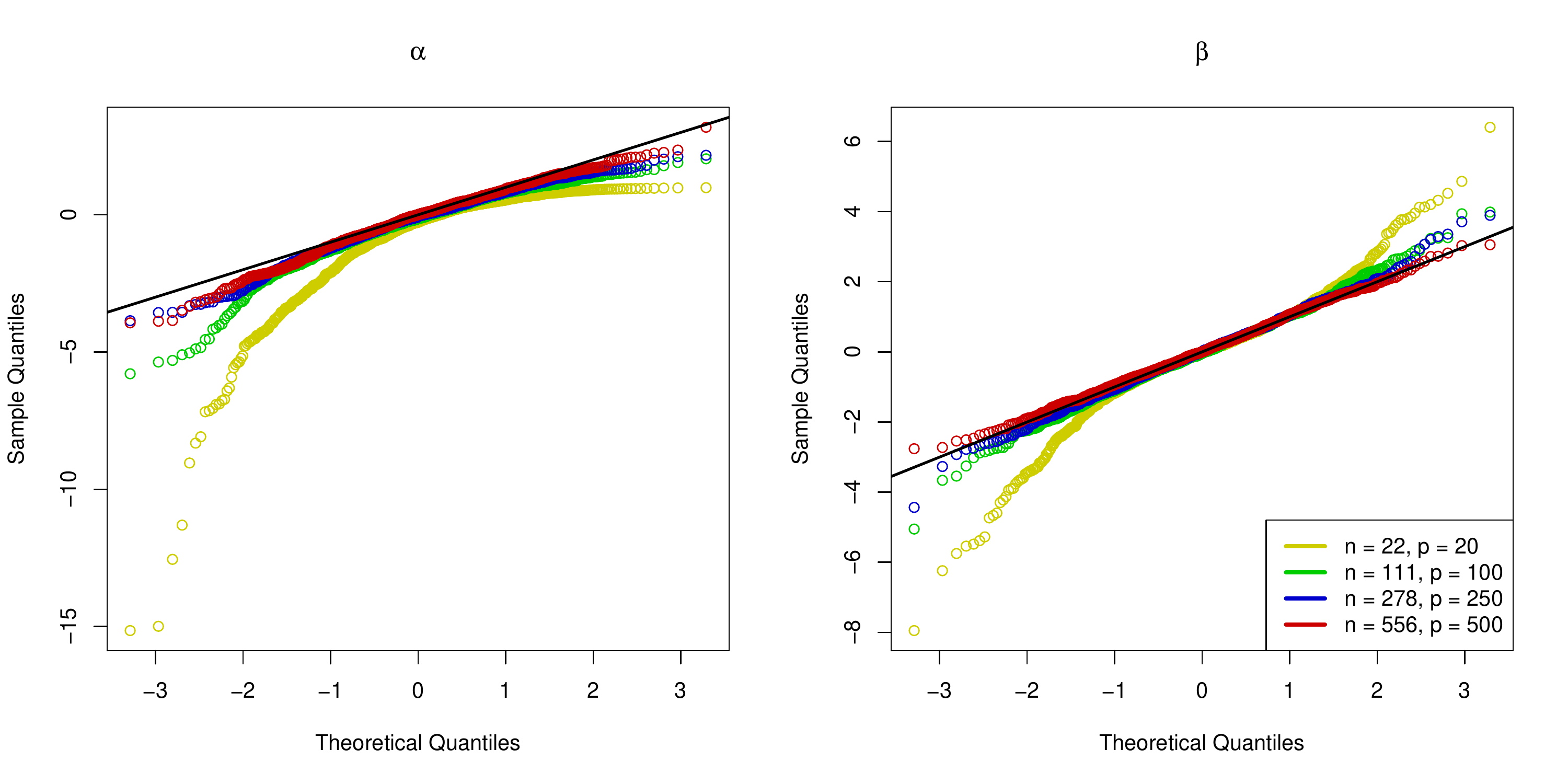}
		\vspace{-5mm}
		\caption{Q-Q plots for $\hat\alpha^*$ and $\hat\beta^*$ in the case of $\gamma = 0$. We set $c = 0.5$ (top panel) and $c = 0.9$ (bottom panel). }\label{fig:norm:bf}
		\end{center}
	\end{figure}

The Q-Q plots for the normal distribution are depicted in Figure~\ref{fig:norm:or}, while Figure~\ref{fig:norm:bf} shows the corresponding results in the case of $\hat\alpha^*$ and $\hat\beta^*$. It is remarkable that Figure~\ref{fig:norm:or} supports the findings of Theorem~\ref{th2}: the constructed Q-Q plots look like straight lines already for $p = 100$ with the exception in the case of $\beta_n^*$ for $c = 2$. Moderate deviations are present only for small dimension $p = 20$ and sample sizes $n = 22$ and $n = 10$ for $c = 0.9$ and $c = 2.0$, respectively. The finite-sample performance of the asymptotic results presented in Theorem~\ref{th4} were investigated in Figure~\ref{fig:norm:bf} for $c<1$. Here, we observed a perfect fit to the normal distribution even for moderate sample sizes and moderate dimensions.

\begin{figure}[!ht]
		\begin{center}
		\includegraphics[scale=0.5]{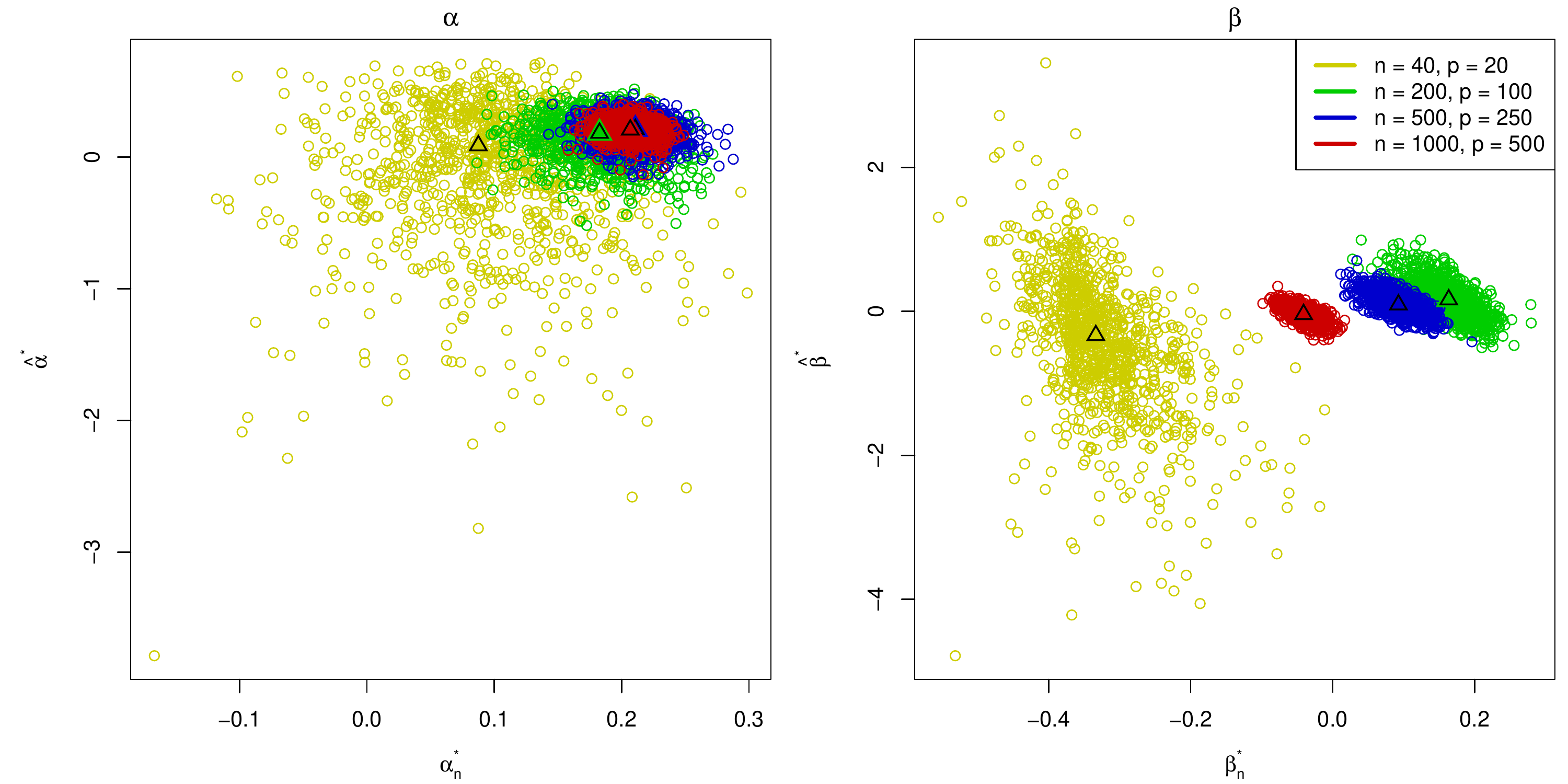}
		\includegraphics[scale=0.5]{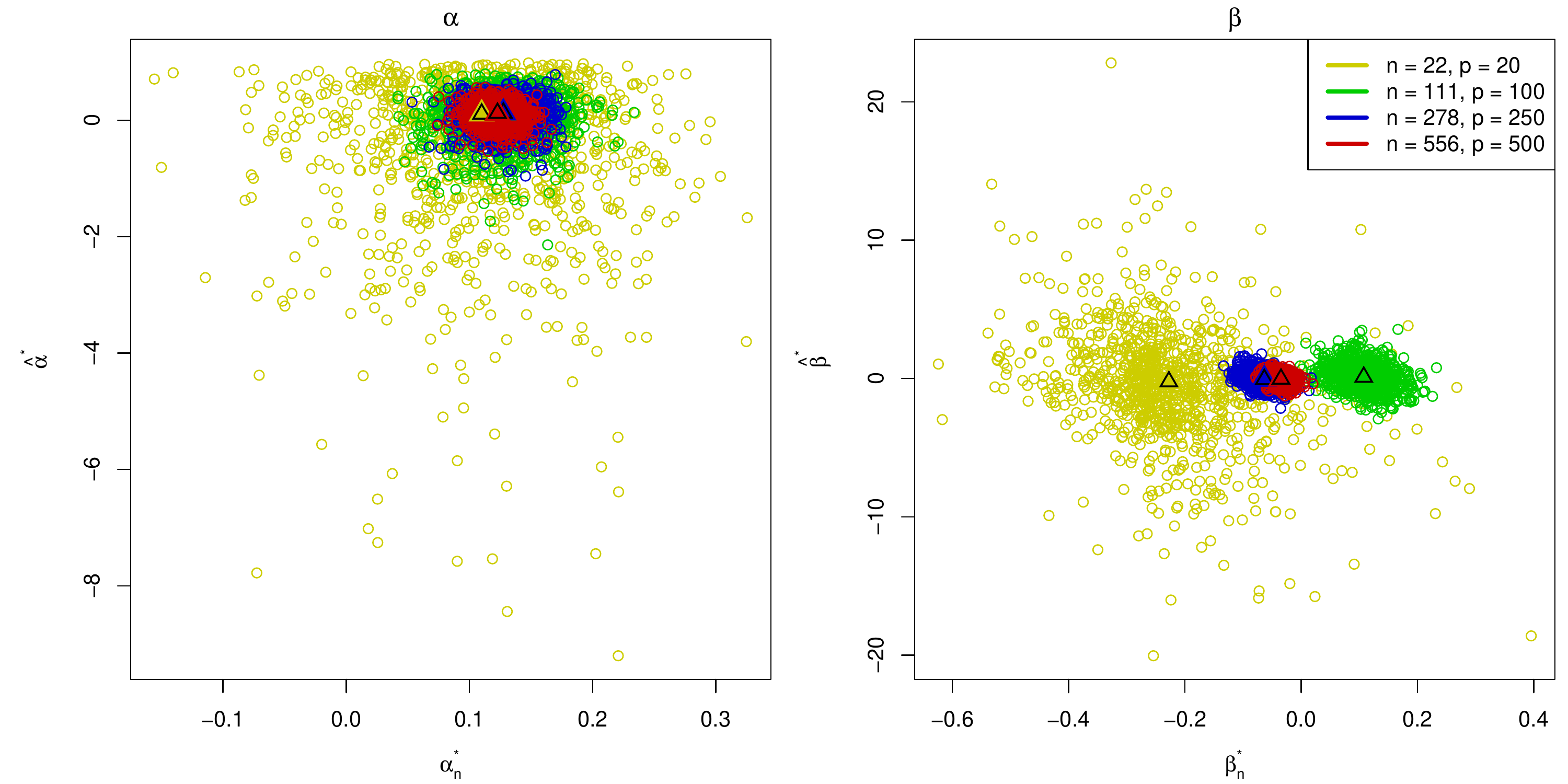}
		\includegraphics[scale=0.5]{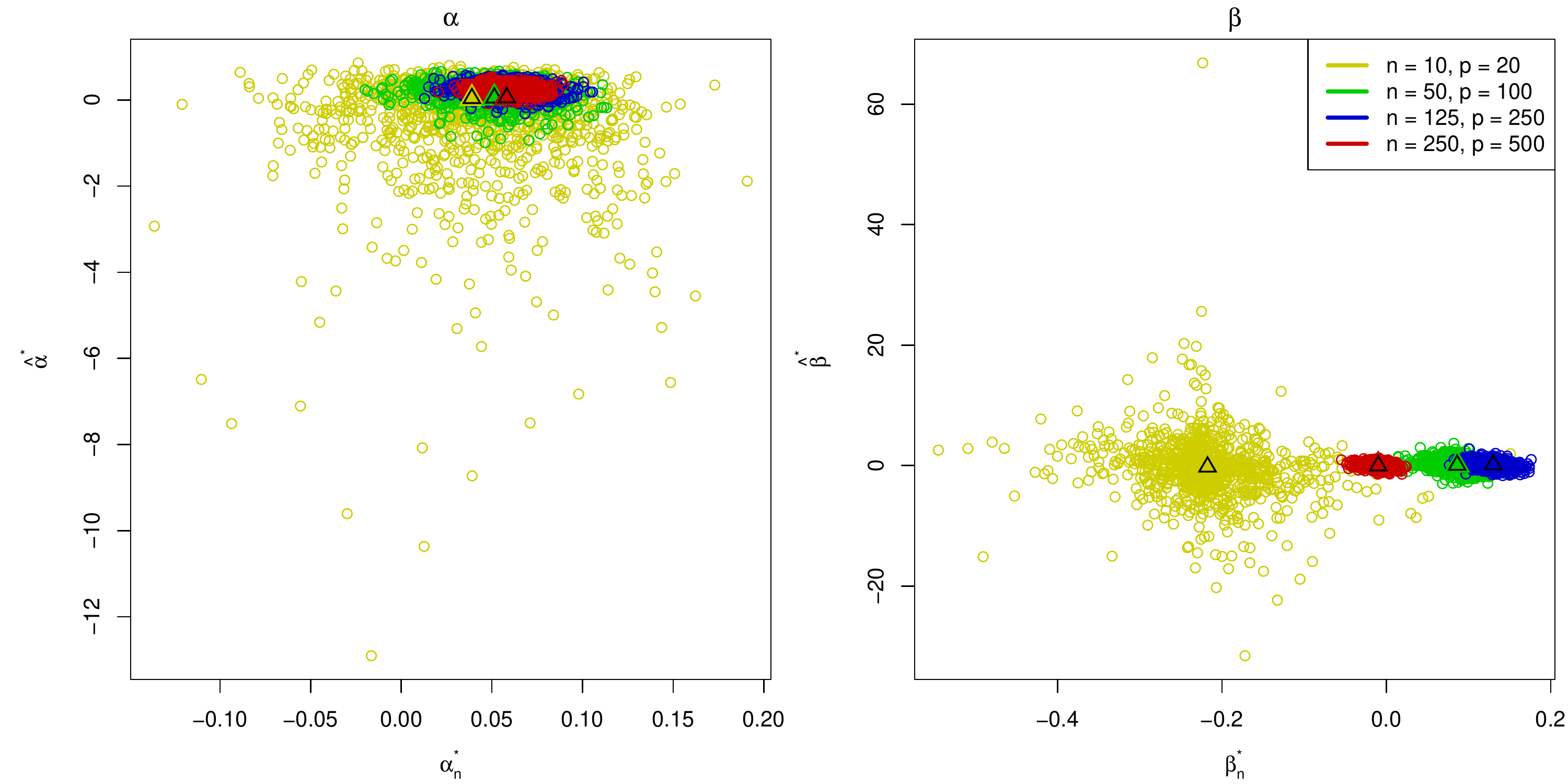}
		\caption{Scatter plots for $(\alpha_n^*, \hat\alpha^*)$ and $(\beta_n^*, \hat\beta^*)$ in the case of $\gamma = 0$. The triangular shows the asymptotic values $(\alpha^*, \alpha^*)$ and $(\beta^*, \beta^*)$. We set $c = 0.5$ (top panel), $c = 0.9$ (middle panel) and $c = 2.0$ (bottom panel).}\label{fig:norm:orVSbf}
		\end{center}
	\end{figure}

The scatter plots of the optimal shrinkage intensities and the corresponding bona fide estimators are provided in Figure~\ref{fig:norm:orVSbf}. We highlighted with triangles the corresponding asymptotic values $\alpha^*$ and $\beta^*$. It is important to note that all the values converge to the asymptotic ones and all the clouds are centered around triangles with the dispersion decreasing with $n$. It is remarkable that $\alpha^*$ should be positive. In contrast, both $\alpha_n$ and $\hat{\alpha}^*$ can be negative. Additional insight is provided by Table~\ref{table:1}, where relative frequencies of $\alpha_n^*$ and $\hat \alpha^*$ being below zero for several values of $p$ and $c$ are presented. It is remarkable that they decrease as $p$ increases.

\begin{figure}[!ht]
		\begin{center}
		\includegraphics[scale=0.7]{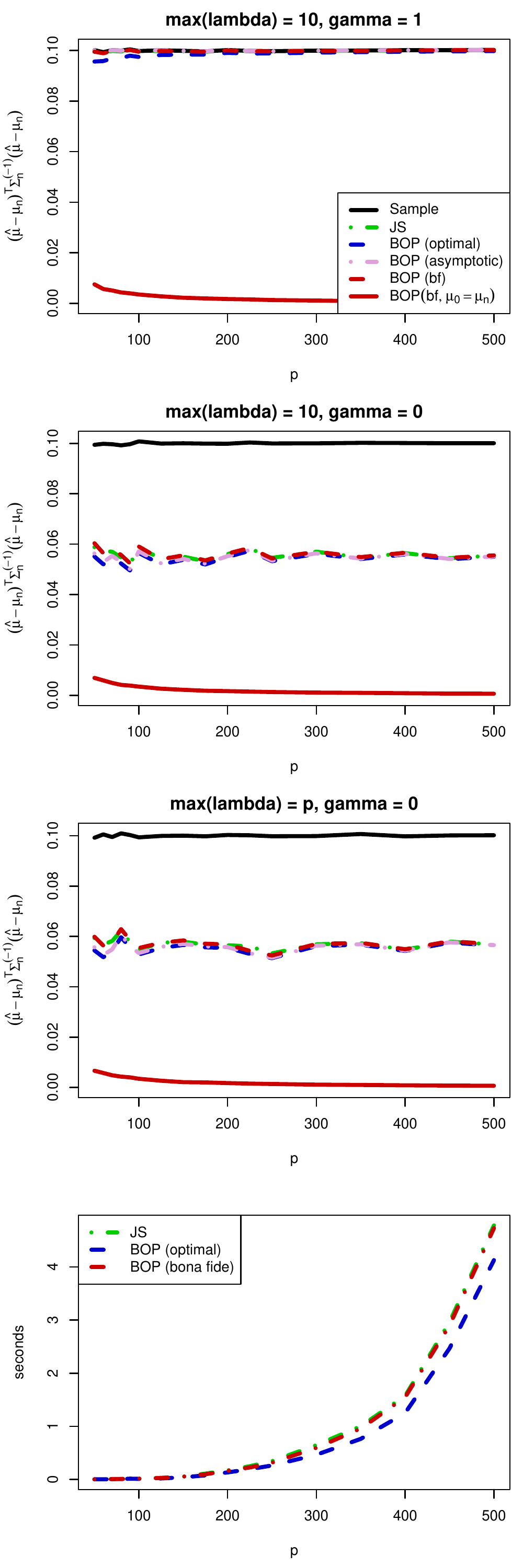}
		\includegraphics[scale=0.7]{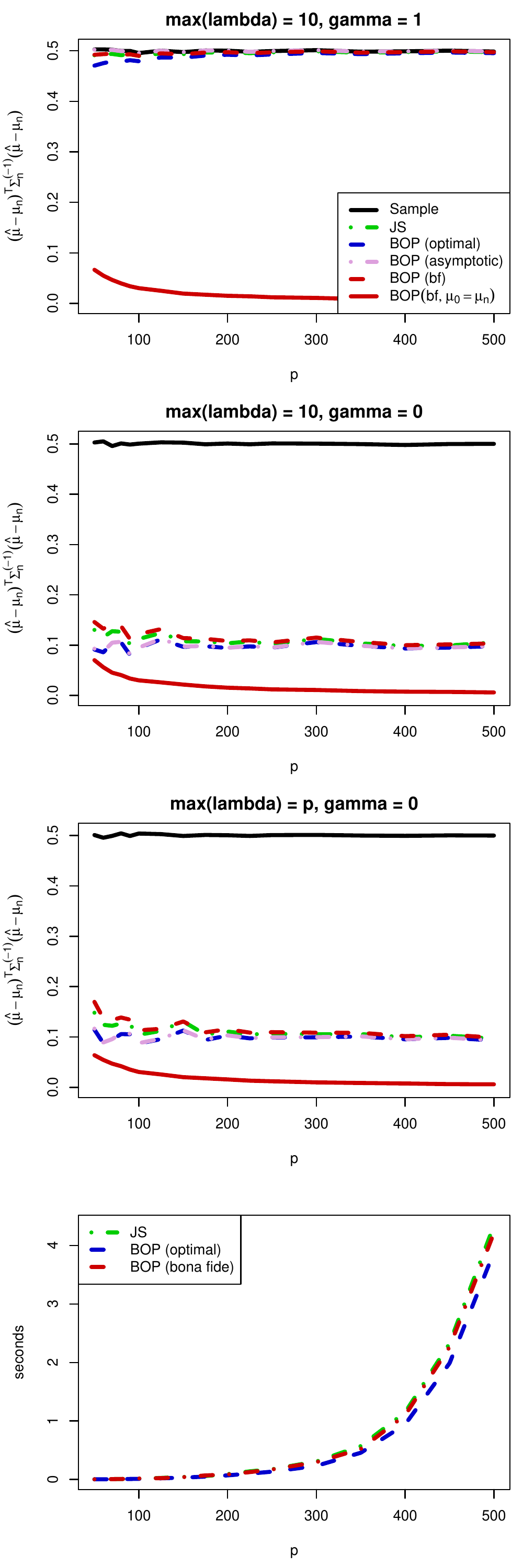}
		\caption{Averaged quadratic loss for estimators $\bby_n$ (black), $\hat{\bm}_n$ (bona fide, red, dashed for $\bm_0\neq\bm_n$ and solid for $\bm_0=\bm_n$), $\bm_n^*$  (optimal, blue), $\bm^*$ (asymptotic, plum), $\widehat{\bm}_{n,JS}$ (green) of $\bm_n$, performed with $N = 1000$ iterations for $c = 0.1$ (left) and $c = 0.5$ (right).}\label{fig:quality:less1}
		\end{center}
	\end{figure}

\begin{figure}[!ht]
		\begin{center}
		\includegraphics[scale=0.7]{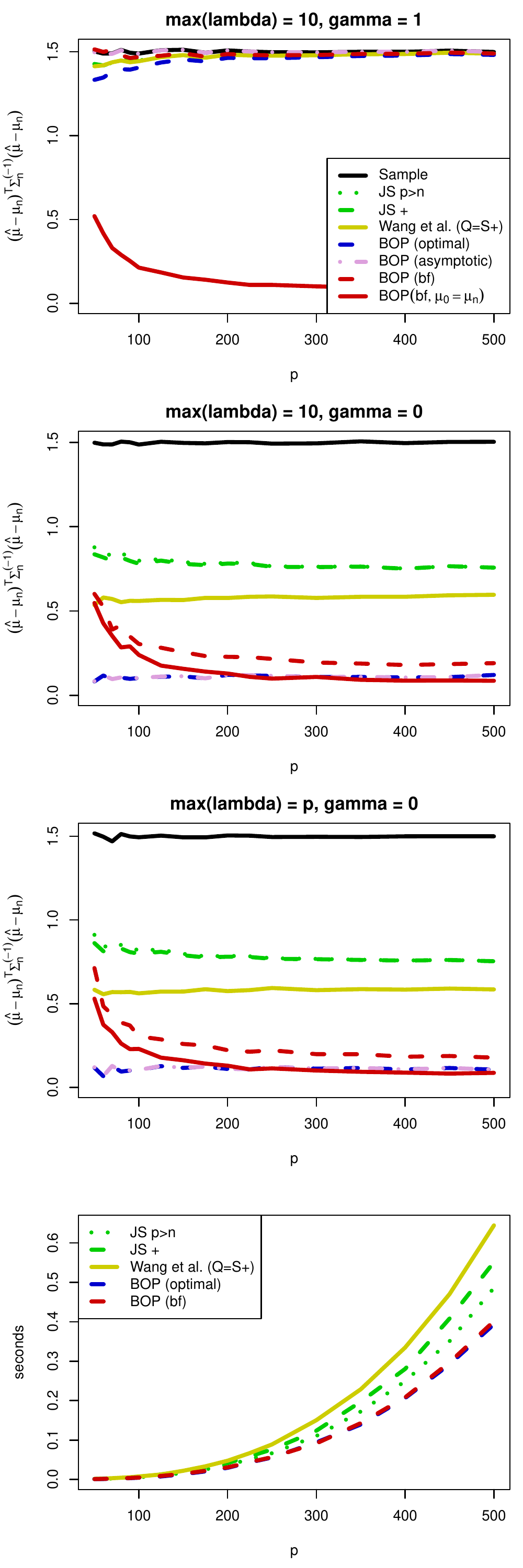}
		\includegraphics[scale=0.7]{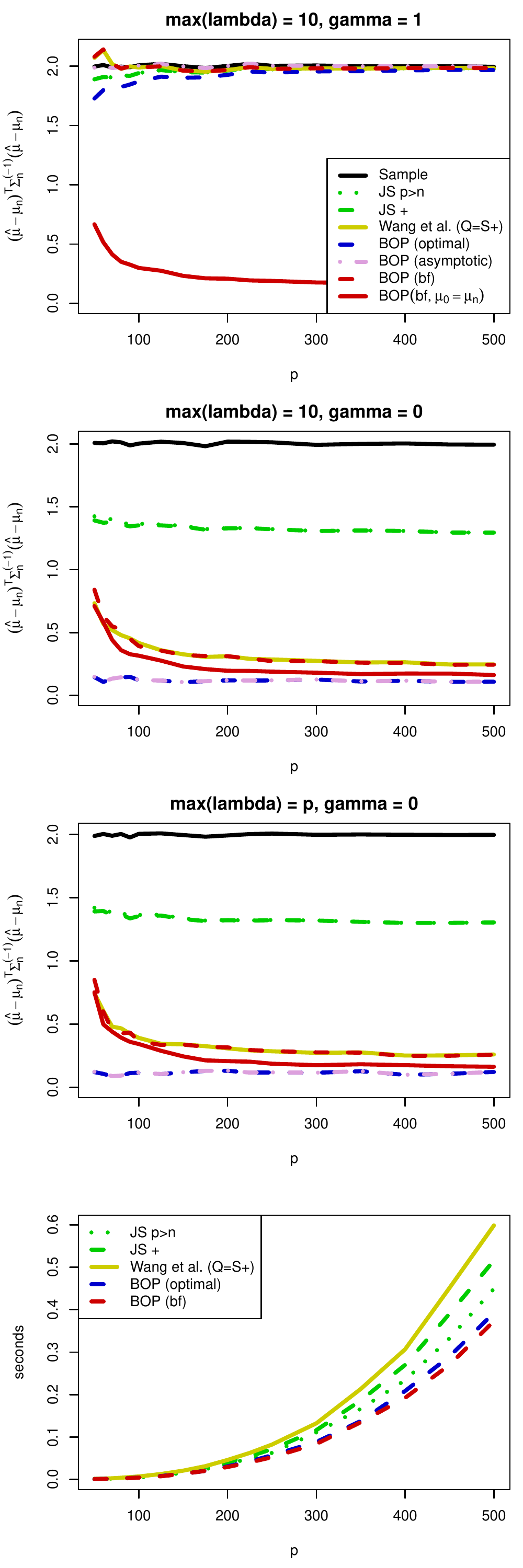}
		\caption{Quadratic loss function for estimators $\bby_n$ (black), $\bm^*$ (asymptotic, plum), $\bm_n^*$  (optimal, blue), $\hat{\bm}_n$  (bona fide, red, dashed for $\bm_0\neq\bm_n$ and solid for $\bm_0=\bm_n$), $\hat{\bm}_{n,W}$ (yellow), $\hat{\bm}_{n,JS(p>n)}$ (green, dotted), $\hat{\bm}_{n,JS+}$ (green, dashed) of $\bm_n$, performed with $N = 1000$ iterations, for $c = 1.5$ (left) and $c = 2.0$ (right).}\label{fig:quality:greater1}
		\end{center}
	\end{figure}

Next, the quality of the estimators was measured by the quadratic loss expressed as
\[
L(\widehat{\bm}, \bm_n, \bSigma_n) = (\widehat{\bm} - \bm_n)^\top\bSigma^{-1}_n(\widehat{\bm} - \bm_n),
\]
where $\widehat{\bm}$ is an estimator for $\bm_n$. The results are depicted in Figures~\ref{fig:quality:less1} and \ref{fig:quality:greater1} in the first, second and  third rows, while the fourth rows present the time in seconds needed to compute each of the estimators (except the asymptotic one). The first and the second rows show the results for the distribution of eigenvalues of $\bSigma_n$ in proportions 20\%, 40\%, 40\% for values 1, 3 and 10 as mentioned in the beginning of this section, while in the third row, $\lambda_{\max}=p$ was chosen which corresponds to the extreme lambda case. The first and second rows differ in the choice of $\gamma$; $\gamma = 1$ in the first row and $\gamma = 0$ in the second row. For each $p$ and $c$ the procedure was repeated $N = 1000$ times and averaged over simulation runs quadratic losses are plotted.

If $\gamma>0$, then $\alpha^* \to  1$ a.s. Hence, the shrinkage estimator tends to the sample mean vector in this case. Interestingly, a similar behavior is also present for all of the considered benchmark estimators. This result is clearly visible in the first rows on both figures. As shown in Figure~\ref{fig:quality:less1}, the optimal shrinkage estimator and its asymptotic counterparts with $\alpha^*_n$ and $\beta^*_n$ replaced by $\alpha^*$ and $\beta^*$ exhibited the best performance. This behavior is not surprising since they both contain unobservable information (true values of $\bSigma_n$ or $\bm_n$). The proposed bona fide estimator appears to be equivalent to the James--Stein estimator for $c < 1$. When $c > 1$, the Moore--Penrose inverse of the sample covariance matrix was employed. The new estimator strongly dominates the rest of the competitors for $c = 1.5$ and shows almost equivalent performance to the $\hat{\bm}_{n,W}$ for $c = 2.0$. It is worth mentioning that the loss function for $\bby_n$ is close to $c$ in all of the considered cases. The proposed estimator is also $30\%$ faster to compute than $\hat{\bm}_{n,W}$. In all cases, we observe very fast convergence of quadratic loss to zero when the target mean vector $\bm_0$ coincides with the true mean vector $\bm_n$, as shown by the red solid lines in Figures \ref{fig:quality:less1} and \ref{fig:quality:greater1}.

\begin{table}[t!]
\caption{Relative frequencies of optimal shrinkage intensity $\alpha*_n$ and its bona fide estimator $\hat \alpha^*$ to be negative.\label{table:1}}

\begin{center}
\begin{tabular}{l|ccc|ccc}\toprule
&                \multicolumn{3}{c|}{$\alpha^*_n$} & \multicolumn{3}{c}{$\hat{\alpha}^*$}\\
$p\setminus c$ &   0.5 &   0.9 &   2.0    &   0.5 &   0.9 &   2.0\\\midrule
20             & 0.082 & 0.071 & 0.180    & 0.505 & 0.555 & 0.621\\
100            & 0.000 & 0.000 & 0.014    & 0.192 & 0.397 & 0.200\\
250            & 0.000 & 0.000 & 0.000    & 0.036 & 0.335 & 0.043\\
500            & 0.000 & 0.000 & 0.000    & 0.006 & 0.268 & 0.010\\\bottomrule
\end{tabular}
\end{center}
\end{table}

\section{Financial application\label{sec:6}}

In order to implement the suggested procedure to real data, we considered weekly log-returns $\mathbf{y}_{t} = (y_{1t}, \ldots, y_{pt})^\top$ of the $p = 412$ constituents of the S\&P500 index traded over $783$ weeks in the period extending from June 12, 2001 to June 9, 2016. We compared the ability of the estimators of the mean vector to predict the future realization of the expected return of the equally-weighted portfolio, which is a popular trading strategy in the high-dimensional case; see, e.g., \cite{demiguel2009optimal}, \cite{tu2011markowitz}. For this reason, we applied the rolling window estimator with window sizes of $n \in\{25, 50, 75, 100\}$ leading to $c \in \{16.48, 8.24, 5.49, 4.12\}$. Namely, for each fixed date $t$ we took historical data of asset returns with the $n$ most recent observations, viz. $\by_{t, n} = \{y_{i\tau}: \tau \in \{t-n, \ldots, t\}, i \in \{1, \ldots, p\}\}$, and estimated the expected return of the equally-weighted portfolio which was then compared to the realized return of this portfolio on day $t+1$.

For the proposed bona fide estimator $\hat\bm_n$(bf), we considered three different choices of the target mean value $\bm_0$. In the first case, the elements of $\bm_0$ is drawn randomly from the uniform distribution on the interval covering the smallest and the largest average of weekly returns for the corresponding stock (over respective rolling window). In the second case, $\bm_0$ consists of randomly interchanging values of $1$ and $-1$, whereas in the third case we set $\bm_0 = \mathbf{1}_p$.

For a given estimator $\hat\bm$ of $\bm_n$, the estimated expected return of the equally-weighted portfolio is calculated by $\hat R_{t, n} (\hat\bm) =  \bi^\top \hat\bm /p$. Then, the performance of each estimator was measured by the average quadratic deviation (times $10^4$ for visualization) from the realized return of the equally-weighted portfolio expressed as
\begin{equation*}
	L = 10^4 \times \frac{1}{T} \sum_{t=1}^T\{\hat R_{t, n} (\hat\bm) - p^{-1}\bi^\top \mathbf{y}_{t+1}\}^2,
\end{equation*}
where $T = 683$ corresponds to the number of rolling windows. The results are depicted in Table~\ref{table:2}. The proposed shrinkage estimator is seen to outperform the other estimators in three out of four cases, particularly in the case where the target mean vector $\bm_0$ contains randomly interchanging values of 1 and $-1$. Only for larger sample size $n = 100$ (with $c = 4.12$), the estimator $\widehat{\bm}_{n,W}$ performs better than the new estimator with any target. In contrast, $\widehat{\bm}_{n,W}$ has the largest value of the average quadratic deviation if $n = 25$. Note also the poor performance of the sample mean vector, which is always one of the worst estimation strategies.

\begin{table}[t!]
\caption{Performance of different mean estimators on the basis of equally weighted portfolio for weekly returns for different sample sizes. The best models are underlined for each sample size.\label{table:2}}

\begin{center}
\begin{tabular}{rr|lrrrrrrrr}
$n$ & $c$ & $\bby_n$ & \multicolumn{2}{c}{$\widehat{\bm}_{JS}$} & $\widehat{\bm}_{n,W}$ & $\widehat{\bm}_n$ (bf) & $\widehat{\bm}_n$ (bf) &$\widehat{\bm}_n$ (bf) \\
& & & $\widehat{\bm}_{JS(p>n)}$ & $\widehat{\bm}_{JS+}$ &  & $\boldsymbol{\mu}_0$ ~  uniform & sample$(1,-1)$  & (1,\ldots,1) &\\
\toprule
$ 25 $ & 16.48 & 7.882 & 7.832 & 7.832 & 7.931 & 7.831 & \underline{7.824} & 7.847 \\
$ 50 $ &  8.24 & 7.813 & 7.721 & 7.721 & 7.748 & 7.716 & \underline{7.715} & 7.736   \\
$ 75 $ &  5.49 & 7.769 & 7.672 & 7.673 & 7.669 & \underline{7.660}& 7.668  & 7.687  \\
$ 100 $ & 4.12 & 7.723 & 7.646 & 7.646 & \underline{7.630} & 7.636 & 7.644 & 7.654 \\
\bottomrule
\end{tabular}
\end{center}
\end{table}

\section{Summary\label{sec:7}}

Nowadays, modern scientific data sets involve sample sizes that are often comparable or even lower than the dimension (number of features) of the data generating mechanism. In this situation many statistical, financial and genetic problems require simple and feasible estimators of the mean vector. Although most of the classical multivariate procedures are based on limit theorems assuming that the dimension $p$ is fixed and the sample size $n$ increases, it has been pointed out by numerous authors that this assumption does not yield precise estimators for commonly used quantities in large dimensions, and that better estimators can be obtained considering scenarios where the dimension tends to infinity as well; see, e.g., \cite{baisil2004} and references therein.

In this paper, we showed how to estimate the mean vector using the optimal shrinkage technique based on random matrix theory. We proved the asymptotic normality of the optimal shrinkage intensities as well as their almost sure limit convergence under the large-dimensional asymptotic regime. Both the simulation study and the empirical application illustrate the dominance of the proposed estimation technique over the benchmark methods in terms of quadratic loss and computational time. Following the recent developments in functional data analysis where a data curve can be treated as a high-dimensional vector of its discretized values on a finite grid \cite{aneiros2015partial, aneiros2016sparse}, the results derived here have potential applications in dealing with estimation problems arising in functional data settings.

\section*{Appendix}

This appendix contains proofs of this paper's results.

\bigskip
\noindent{\it Proof of Theorem~\ref{th1}}.
 Let
\begin{align}
\eta_1&=p^{-\gamma}\bby_n^\top\bSigma^{-1}_n\bxi=n^{-1}p^{-\gamma}\bi^\top_n\by_n^\top\bSigma^{-1}_n\bxi
=n^{-1}p^{-\gamma}\bi^\top_n\bx_n^\top\bSigma_n^{-1/2}\bxi+p^{-\gamma}\bm_n^\top\bSigma_n^{-1}\bxi,\label{app_th1_eq1}\\
\eta_2&=p^{-\gamma}\bby_n^\top\bSigma^{-1}_n\bby_n=n^{-2}p^{-\gamma} \bi^\top_n \by^\top_n\bSigma^{-1}_n\by_n\bi_n \nonumber \\
&=n^{-2}p^{-\gamma}\bi^\top_n(\bx_n^\top\bx_n)\bi_n+2n^{-1}p^{-\gamma}\bi^\top_n\bx_n^\top\bSigma_n^{-1/2}\bm_n+p^{-\gamma}\bm_n^\top\bSigma_n^{-1}\bm_n,\label{app_th1_eq2}
\end{align}
where the vector $\bxi$ is either $\bm_0$ or $\bm_n$. Then
\begin{equation}\label{app_th1_eq3}
n^{-1}p^{-\gamma}\bi^\top_n\bx_n^\top\bSigma_n^{-1/2}\bxi = n^{-1}p^{-\gamma} (\bi_n^\top \otimes \bxi^\top \bSigma_n^{-1/2}){\rm vec} (\bx_n)
\end{equation}
and
$n^{-2} \bi^\top_n(\bx_n^\top\bx_n)\bi_n =  \|  \bx_n \bi_n /\sqrt{n} \|^2 /n$.

In Eq.~\eqref{app_th1_eq3}, the sum of $pn$ independent random variables is present with zero means and variances
equal to $(\bxi^\top\bSigma_n^{-1/2}\be_i)^2$ for $i \in \{ 1,\ldots , p\}$, where $\be_i$ is the $i$th basis vector in $\mathbb{R}^p$. Because
\begin{align*}
{\rm var} \{n^{-1}p^{-\gamma} (\bi_n^\top \otimes \bxi^\top \bSigma_n^{-1/2}){\rm vec} (\bx_n) \} &= \displaystyle \frac{p^{-2\gamma}n}{n^2} \sum_{i=1}^p (\bxi^\top\bSigma_n^{-1/2}\be_i)^2= \frac{p^{-2\gamma}}{n} \bxi^\top\bSigma_n^{-1}\bxi\\
& \displaystyle \le \lambda_{\max}(\bSigma_n^{-1})\frac{p^{-\gamma}}{n} p^{-\gamma}\|\bxi\|^2 \le \frac{\lambda_{\max}(\bSigma_n^{-1})M_u}{np^{\gamma}} < \infty,
\end{align*}
the application of Kolmogorov's Strong Law of Large Numbers implies that, as $p,n \rightarrow\infty$,
$p^{-\gamma}n^{-1}\bi^\top_n\bx_n^\top\bSigma_n^{-1/2}\bxi \stackrel{\rm a.s.}{\to }0$.

Next, we consider the elements of the vector $\bx_n\bi_n$ which are determined by the sums of
independently and identically distributed random variables with zero means and unit variances.
Hence, from the Central Limit Theorem the components of the vector $n^{-1/2}\bx_n\bi_n$ are asymptotically independent
and standard Gaussian. Consequently, their squares are asymptotically independent and $\chi^2$-
distributed with one degree of freedom. Because the expectation of a $\chi^2_1$-distributed random
variable is $1$ and the variance is finite we get that, as $p, n \to \infty$,
$\|  \bx_n \bi_n /\sqrt{n} \|^2 /p \stackrel{\rm a.s.}{\to }1$
and hence, as $p, n \to \infty$,
$n^{-2}\bi^\top_n(\bx_n^\top\bx_n)\bi_n \stackrel{\rm a.s.}{\to }c$.
Hence,
\begin{align*}
\alpha_n^*&= \frac{\bby_n^\top\bSigma^{-1}_n\bm_n\bm_0^\top\bSigma^{-1}_n\bm_0-\bm_n^\top\bSigma^{-1}_n\bm_0\bby_n^\top\bSigma^{-1}_n\bm_0}
{\bby_n^\top\bSigma^{-1}_n\bby_n\bm_0^\top\bSigma^{-1}_n\bm_0-(\bby_n^\top\bSigma^{-1}_n\bm_0)^2}
= \frac{p^{-\gamma}\bby_n^\top\bSigma^{-1}_n\bm_n p^{-\gamma} \bm_0^\top\bSigma^{-1}_n\bm_0- p^{-\gamma} \bm_n^\top\bSigma^{-1}_n\bm_0 p^{-\gamma} \bby_n^\top\bSigma^{-1}_n\bm_0}
{p^{-\gamma} \bby_n^\top\bSigma^{-1}_n\bby_n p^{-\gamma} \bm_0^\top\bSigma^{-1}_n\bm_0-(p^{-\gamma}\bby_n^\top\bSigma^{-1}_n\bm_0)^2}\\
&\rightarrow \frac{p^{-\gamma}\bm_n^\top\bSigma^{-1}_n\bm_n p^{-\gamma}\bm_0^\top\bSigma^{-1}_n\bm_0- p^{-\gamma} \bm_n^\top\bSigma^{-1}_n\bm_0 p^{-\gamma}\bm_n^\top\bSigma^{-1}_n\bm_0}
{(p^{-\gamma}c+ p^{-\gamma}\bm_n^\top\bSigma^{-1}_n\bm_n) p^{-\gamma} \bm_0^\top\bSigma^{-1}_n\bm_0-(p^{-\gamma}\bm_n^\top\bSigma^{-1}_n\bm_0)^2} =  \frac{\bm_n^\top\bSigma^{-1}_n\bm_n \bm_0^\top\bSigma^{-1}_n\bm_0-  (\bm_n^\top\bSigma^{-1}_n\bm_0)^2}
{(c+ \bm_n^\top\bSigma^{-1}_n\bm_n) \bm_0^\top\bSigma^{-1}_n\bm_0-(\bm_n^\top\bSigma^{-1}_n\bm_0)^2}=\alpha^*,
\end{align*}
and, similarly,
\begin{align*}
\beta_n^*&= \frac{\bby_n^\top\bSigma^{-1}_n\bby_n\bm_n^\top\bSigma^{-1}_n\bm_0-\bby_n^\top\bSigma^{-1}_n\bm_0\bby_n^\top\bSigma^{-1}_n\bm_n}
{\bby_n^\top\bSigma^{-1}_n\bby_n\bm_0^\top\bSigma^{-1}_n\bm_0-(\bby_n^\top\bSigma^{-1}_n\bm_0)^2}
= \frac{p^{-\gamma}\bby_n^\top\bSigma^{-1}_n\bby_n p^{-\gamma}\bm_n^\top\bSigma^{-1}_n\bm_0-p^{-\gamma} \bby_n^\top\bSigma^{-1}_n \bm_0 p^{-\gamma}\bby_n^\top\bSigma^{-1}_n\bm_n}
{p^{-\gamma}\bby_n^\top\bSigma^{-1}_n\bby_n p^{-\gamma} \bm_0^\top\bSigma^{-1}_n\bm_0-(p^{-\gamma}\bby_n^\top\bSigma^{-1}_n\bm_0)^2}\\
&\rightarrow \frac{(c+ \bm_n^\top\bSigma^{-1}_n\bm_n)\bm_n^\top\bSigma^{-1}_n\bm_0-  \bm_n^\top\bSigma^{-1}_n\bm_0 \bm_n^\top\bSigma^{-1}_n\bm_n}
{(c+ \bm_n^\top\bSigma^{-1}_n\bm_n)  \bm_0^\top\bSigma^{-1}_n\bm_0-(\bm_n^\top\bSigma^{-1}_n\bm_0)^2}
=(1-\alpha^*)\frac{\bm_n^\top\bSigma^{-1}_n\bm_0}{\bm_0^\top\bSigma^{-1}_n\bm_0}.
\end{align*}
for $p/n \rightarrow c>0$ as $n \rightarrow\infty$. The theorem is thus proved. \hfill $\Box$

\bigskip
In the proof of Theorem~\ref{th2} we make use of the following lemma.

\begin{lemma}\label{lem1}
Let $\bZ \sim \mathcal{N}_m(\mathbf{0},\bI_m)$ and define
$Q_1(\bZ)=\sqrt{m} \, ( \bZ^\top/m \bZ-1 )$ and $Q_2(\bZ)= \mathbf{a}^\top \bZ/\sqrt{m}$ for nonrandom $ \ba$ with $\lim_{m\rightarrow \infty} {\ba^\top\ba}/{m}=\tilde{\sigma}^2$.
Then, as $m\rightarrow\infty$,
$$
\left(
  \begin{array}{c}
    Q_1(\bZ) \\
    Q_2(\bZ) \\
  \end{array}
\right) \rightsquigarrow \mathcal{N}_2\left[\mathbf{0},\left(
                                                      \begin{array}{cc}
                                                        2 & 0 \\
                                                        0 & \tilde{\sigma}^2\\
                                                      \end{array}
                                                    \right) \right].
$$
\end{lemma}

\noindent{\it Proof of Lemma~\ref{lem1}.}
From Lemma 2.1 of \cite{schoeneschmid2000} we get the moment generating function of $(Q_1(\bZ),Q_2(\bZ))^\top$ expressed as
\begin{align*}
M(t_1,t_2)= \operatorname{E} [\exp\{t_1Q_1(\bZ)+t_2Q_2(\bZ)\} ] = \exp(-\sqrt{m}t_1) (1- {2t_1}/{\sqrt{m}} )^{-m/2}\exp\left(\frac{t_2^2}{2}\frac{\ba^\top\ba/m}{1-2t_1/\sqrt{m}}\right).
\end{align*}

Hence,
\begin{align*}
\lim_{m\rightarrow \infty} M(t_1,t_2)&=  \exp ( {\tilde{\sigma}^2 t_2^2}/{2} ) \lim_{m\rightarrow \infty}
\exp(-\sqrt{m}t_1)  (1- {2t_1}/{\sqrt{m}} )^{-m/2}=\exp ({\tilde{\sigma}^2 t_1^2}/{2}+t_1^2 )
,
\end{align*}
where the last equality follows from the facts that the function under the limit is the moment generating function of $\sqrt{m}\, (\chi/m-1)$ with $\chi\sim\chi^2_m$ and $\sqrt{m} \, (\chi/m-1) \rightsquigarrow  \mathcal{N}(0,2)$ as $m\rightarrow\infty$. The lemma is thus proved. \hfill $\Box$

\bigskip
\noindent{\it Proof of Theorem~\ref{th2}.}
Let $\tilde{c}= p^{-\gamma}c$, $q_{ij} = p^{-\gamma}\bm_i^\top\bSigma^{-1}_n\bm_j$, for $i, j\in\{0, n\}$ and $d = q_{00}q_{nn}-q_{0n}^2$.
First, we prove the result in the case of $\alpha_n^*$. Then
\begin{align*}
\sqrt{p^{\gamma}n}\, (\alpha_n^*-\alpha^*)&=
\frac{1}{(\tilde{c}\qnullnull +d)^2}\frac{\tilde{c}\qnullnull +d}{\qnullnull p^{-\gamma}\bby_n^\top\bSigma^{-1}_n\bby_n-(p^{-\gamma}\bm_0^\top\bSigma^{-1}_n\bby_n)^2}Z_{\alpha},
\end{align*}
where
\begin{align*}
Z_{\alpha}= \sqrt{p^{\gamma}n} \, [(\tilde{c}\qnullnull +d)(\qnullnull p^{-\gamma}\bm_n^\top\bSigma^{-1}_n\bby_n-\qnulln \bm_0^\top\bSigma^{-1}_n\bby_n) - d\{\qnullnull p^{-\gamma}\bby_n^\top\bSigma^{-1}_n\bby_n-(p^{-\gamma}\bm_0^\top\bSigma^{-1}_n\bby_n)^2\} ].
\end{align*}
The application of \eqref{app_th1_eq1} and \eqref{app_th1_eq2} leads to
\begin{align*}
Z_{\alpha}&= \sqrt{p^{\gamma}n}\big((\tilde{c}\qnullnull +d)\{\qnullnull (\qnn +n^{-1}p^{-\gamma}\bm_n^\top\bSigma^{-1/2}_n\bx_n \bi_n)
-\qnulln (\qnulln +n^{-1}p^{-\gamma}\bm_0^\top\bSigma^{-1/2}_n\bx_n\bi_n)\}\\
& \quad \quad - d[\qnullnull \{\qnn +2n^{-1}p^{-\gamma}\bm_n^\top\bSigma_n^{-1/2}\bx_n\bi_n +n^{-2}p^{-\gamma}\bi^\top_n(\bx_n^\top\bx_n)\bi_n\}\\
& \quad \quad -\{\qnulln ^2+2\qnulln n^{-1}p^{-\gamma}\bm_0^\top\bSigma^{-1/2}_n\bx_n\bi_n+n^{-2}p^{-2\gamma}(\bm_0^\top\bSigma^{-1/2}_n\bx_n\bi_n)^2 \}]\big)\\
&= [(\tilde{c}\qnullnull -d)(\qnullnull p^{-\gamma/2}\bm_n^\top\bSigma^{-1/2}_n-\qnulln p^{-\gamma/2}\bm_0^\top\bSigma^{-1/2}_n)(n^{-1/2}\bx_n\bi_n) .\\
&\quad \quad +dp^{-\gamma/2}n^{-1/2}\{p^{-\gamma/2}\bm_0^\top\bSigma^{-1/2}_n (n^{-1/2}\bx_n\bi_n)\}^2
- \sqrt{\tilde{c}}\, d\qnullnull \sqrt{p} \, ( \|n^{-1/2}\bx_n\bi_n\|^2/p-1 )].
\end{align*}
Because the elements of $n^{-1/2}\bx_n\bi_n$ are independent and standard Gaussian, we get that, as $p, n \to \infty$,
$$
p^{-\gamma/2}\bm_0^\top\bSigma^{-1/2}_n  (n^{-1/2}\bx_n\bi_n )\rightsquigarrow  \mathcal{N}(0,\qnullnull )\rightarrow\infty
$$
and, consequently,  as $p, n \to \infty$,
\[
\frac{1}{\sqrt{p^{\gamma}n}} \{p^{-\gamma/2}\bm_0^\top\bSigma^{-1/2}_n (n^{-1/2}\bx_n\bi_n) \}^2
\stackrel{\rm a.s.}{\longrightarrow} 0.
\]
Furthermore, the application of Lemma~\ref{lem1} leads to
$Z_{\alpha} \rightsquigarrow  \mathcal{N} (0,\sigma^2_{Z_{\alpha}} )$ for $ p/n\to  c>0$ as $n\to \infty$
with
\begin{align*}
\sigma^2_{Z_{\alpha}}&=  (\tilde{c}\qnullnull -d)^2(\qnullnull ^2 \qnn +\qnulln ^2\qnullnull -2\qnullnull \qnulln ^2)+\tilde{c}d^2\qnullnull ^2
= (\tilde{c}\qnullnull -d)^2\qnullnull d+\tilde{c}d^2\qnullnull ^2.
\end{align*}
Finally, from the proof of Theorem~\ref{th1} we get that, for $p/n\to  c>0$ as $n\to \infty$,
\[\frac{\tilde{c}\qnullnull +d}{\qnullnull p^{-\gamma}\bby_n^\top\bSigma^{-1}_n\bby_n-(p^{-\gamma}\bm_0^\top\bSigma^{-1}_n\bby_n)^2}
\stackrel{\rm a.s.}{\longrightarrow} 1,
\]
which together with Slutsky's theorem (see \cite{lehmann1999}) implies that, for $p/n\to  c>0$ as $n\to \infty$,
\begin{equation*}
\sqrt{p^{\gamma}n} \, (\alpha_n^*-\alpha^*)\rightsquigarrow  \mathcal{N}\left[0,\frac{(\tilde{c}\qnullnull -d)^2\qnullnull d+\tilde{c}d^2\qnullnull ^2}{(\tilde{c}\qnullnull +d)^4}\right].
\end{equation*}

Similarly, for $\beta_n^*$ we get
\begin{align*}
\sqrt{p^{\gamma}n} \, (\beta_n^*-\beta^*)&=
\frac{1}{(\tilde{c}\qnullnull +d)^2}\frac{\tilde{c}\qnullnull +d}{\qnullnull p^{-\gamma}\bby_n^\top\bSigma^{-1}_n\bby_n-(p^{-\gamma}\bm_0^\top\bSigma^{-1}_n\bby_n)^2}Z_{\beta}
\end{align*}
with
\begin{align*}
Z_{\beta}&=
\sqrt{p^{\gamma}n}\left[(\tilde{c}\qnullnull +d)(p^{-\gamma}\bby_n^\top\bSigma^{-1}_n\bby_n \qnulln -p^{-\gamma}\bby_n^\top\bSigma^{-1}_n\bm_0
p^{-\gamma}\bby_n^\top\bSigma^{-1}_n\bm_n)\right. -\left.\tilde{c }\qnulln \{p^{-\gamma}\bby_n^\top\bSigma^{-1}_n\bby_n\qnullnull -(p^{-\gamma}\bby_n^\top\bSigma^{-1}_n\bm_0)^2\}\right]\\
&= \sqrt{p^{\gamma}n}\{d\qnulln p^{-\gamma}\bby_n^\top\bSigma^{-1}_n\bby_n+\tilde{c} \qnulln (p^{-\gamma}\bby_n^\top\bSigma^{-1}_n\bm_0)^2
-(\tilde{c}\qnullnull +d)p^{-\gamma}\bby_n^\top\bSigma^{-1}_n\bm_0
p^{-\gamma}\bby_n^\top\bSigma^{-1}_n\bm_n\}.
\end{align*}
Using (\ref{app_th1_eq1}) and (\ref{app_th1_eq2}) we get
\begin{align*}
Z_{\beta}&= \sqrt{p^{\gamma}n}[
d\qnulln \{\qnn +2n^{-1}p^{-\gamma}\bm_n^\top\bSigma_n^{-1/2}\bx_n\bi_n +n^{-2}p^{-\gamma}\bi^\top_n(\bx_n^\top\bx_n)\bi_n\}\\
&\quad \quad  +\tilde{c} \qnulln \{\qnulln ^2+2\qnulln n^{-1}p^{-\gamma}\bm_0^\top\bSigma^{-1/2}_n\bx_n\bi_n+n^{-2}p^{-2\gamma}(\bm_0^\top\bSigma^{-1/2}_n\bx_n\bi_n)^2 \}\\
&\quad \quad - (\tilde{c}\qnullnull +d)\{\qnulln \qnn +\qnulln n^{-1}p^{-\gamma}\bm_n^\top\bSigma^{-1/2}_n\bx_n\bi_n\\
& \quad \quad  +\qnn n^{-1}p^{-\gamma}\bm_0^\top\bSigma^{-1/2}_n\bx_n\bi_n
+n^{-2}p^{-2\gamma}(\bm_0^\top\bSigma^{-1/2}_n\bx_n\bi_n)(\bm_n^\top\bSigma^{-1/2}_n\bx_n\bi_n)\}]\\
&= \{(d-\tilde{c}\qnullnull )\qnulln p^{-\gamma/2}\bm_n^\top\bSigma^{-1/2}_n
+(\tilde{c}\qnulln ^2-\tilde{c}d-d\qnn )p^{-\gamma/2}\bm_0^\top\bSigma^{-1/2}_n\} (n^{-1/2}\bx_n\bi_n )\\
& \quad \quad  +\sqrt{\tilde{c}}d\qnulln \sqrt{p} ( \|n^{-1/2}\bx_n\bi_n\|^2/p-1 )
+\tilde{c} \qnulln  n^{-2}p^{-2\gamma}(\bm_0^\top\bSigma^{-1/2}_n\bx_n\bi_n)^2\\
&\quad \quad  -(\tilde{c}\qnullnull +d)n^{-2}p^{-2\gamma}(\bm_0^\top\bSigma^{-1/2}_n\bx_n\bi_n)(\bm_0^\top\bSigma^{-1/2}_n\bx_n\bi_n)
.
\end{align*}
Proceeding as in the proof pertaining to $\alpha_n^*$, we get
\[
\tilde{c} \qnulln \, \frac{p^{-2\gamma}}{n^2}\, (\bm_0^\top\bSigma^{-1/2}_n\bx_n\bi_n)^2
-(\tilde{c}\qnullnull +d)\, \frac{p^{-2\gamma}}{n^2} \, (\bm_0^\top\bSigma^{-1/2}_n\bx_n\bi_n)(\bm_0^\top\bSigma^{-1/2}_n\bx_n\bi_n)\stackrel{\rm a.s.}{\longrightarrow}0
\]
for $p/n\to  c>0$ as $n\to \infty$. Then, the application of Lemma~\ref{lem1} leads to
$Z_{\beta} \rightsquigarrow  \mathcal{N} (0,\sigma^2_{Z_{\beta}} )$ for $p/n\to  c>0$ as $n\to \infty$
with
\begin{align*}
\sigma^2_{Z_{\beta}}&=  (d-\tilde{c}\qnullnull )^2\qnulln ^2\qnn +(\tilde{c}\qnulln ^2-\tilde{c}d-d\qnn )^2\qnullnull
+2(\tilde{c}\qnulln ^2-\tilde{c}d-d\qnn )(d-\tilde{c}\qnullnull )\qnulln ^2+\tilde{c}d^2\qnulln ^2.
\end{align*}
Finally, using Slutsky's theorem (see \cite{lehmann1999}) we get
$\sqrt{p^{\gamma}n}\,(\beta_n^*-\beta^*)\rightsquigarrow  \mathcal{N} [ 0, {\sigma^2_{Z_{\beta}}}/{(\tilde{c}\qnullnull +d)^4} ]$
for $p/n\to  c>0$ as $n\to \infty$. The theorem is proved. \hfill $\Box$

\bigskip
In the proof of Theorem~\ref{th3} we make use of the following lemma.

\begin{lemma}\label{lem2}
Assume (A1)--(A2) and let the elements of $\bx_n$ possess uniformly bounded $4+\varepsilon$, $\varepsilon>0$, moments. Let $\boldsymbol{\theta}$ and $\boldsymbol{\xi}$ be the universal nonrandom vectors from the set $  \{\bm_0, \bm_n \}$. Then Then
$$
   p^{-\gamma}\left|\boldsymbol{\xi}^\top\bS_n^{-1}\boldsymbol{\theta}- (1-c)^{-1} \boldsymbol{\xi}^\top \bSigma^{-1}_n\boldsymbol{\theta}\right| \stackrel{\rm a.s.}{\longrightarrow} 0 ,\quad
   \bbx_n^\top\bSigma_n^{1/2}\bS_n^{-1}\bSigma_n^{1/2}\bbx_n \stackrel{\rm a.s.}{\longrightarrow} \frac{c}{1-c}  ,\quad
 p^{-\gamma/2}\bbx_n^\top\bSigma_n^{1/2}\bS_n^{-1}\boldsymbol{\theta} \stackrel{\rm a.s.}{\longrightarrow} 0
$$
for $p/n\to  c \in (0, \infty)$ as $n\rightarrow\infty$, where $\bbx_n=n^{-1}\bx_n \bi_n$ stands for the sample mean vector calculated from $\bx_n$.
\end{lemma}

\noindent{\it Proof of Lemma~\ref{lem2}.}
For $\boldsymbol{\theta} \in \{\bm_0, \bm_n \}$, one has
$\|\bSigma^{-1/2}_n\boldsymbol{\theta} \|^2\le \lambda_{\max}(\bSigma^{-1}_n) \|\boldsymbol{\theta} \|^2<\infty$
in view of assumptions (A1) and (A2). The rest of the proof follows from Lemma~5.3 in \cite{bodnar2016}. \hfill $\Box$

\bigskip
\noindent{\it Proof of Theorem~\ref{th3}.}
For $\boldsymbol{\theta},\bxi \in \{\bm_0, \bm_n \}$, the application of Lemma~\ref{lem2} leads to
$$
p^{-\gamma} |\bm_0^\top\bS_n^{-1}\bm_0- (1-c)^{-1} \bm_0^\top \bSigma^{-1}_n\bm_0 t| \stackrel{\rm a.s.}{\longrightarrow} 0,
$$
$$
p^{-\gamma} |\bm_0^\top\bS_n^{-1}\bby_n-(1-c)^{-1} \bm_0^\top \bSigma^{-1}_n\bm_n | \le p^{-\gamma} |\bm_0^\top\bS_n^{-1}\bSigma_n^{1/2}\bbx_n |
+ p^{-\gamma} |\bm_0^\top\bS_n^{-1}\bm_n-(1-c)^{-1}\bm_0^\top \bSigma^{-1}_n\bm_n |  \stackrel{\rm a.s.}{\longrightarrow} 0,
$$
and
\begin{multline*}
p^{-\gamma} |\bby_n^\top\bS_n^{-1}\bby_n-(1-c)^{-1} \bm_n^\top \bSigma^{-1}_n\bm_n-c (1-c)^{-1} | \le p^{-\gamma}
|\bm_n^\top\bS_n^{-1}\bm_n-(1-c)^{-1} \bm_n^\top \bSigma^{-1}_n\bm_n |\\ +p^{-\gamma} |\bbx_n^\top\bSigma_n^{1/2}\bS_n^{-1}\bSigma_n^{1/2}\bbx_n-c(1-c)^{-1} |
+2p^{-\gamma} |\bm_n^\top\bS_n^{-1}\bSigma_n^{1/2}\bbx_n |\stackrel{\rm a.s.}{\longrightarrow} 0.
\end{multline*}
Hence,
$$
\hat{\alpha}^* = \frac{\{\bby_n^\top \bS^{-1}_n\bby_n -c/(1-c) \}\bm_0^\top \bS^{-1}_n\bm_0-(\bby_n^\top \bS^{-1}_n\bm_0)^2}{\bby_n^\top \bS^{-1}_n\bby_n \bm_0^\top \bS^{-1}_n\bm_0-(\bby_n^\top \bS^{-1}_n\bm_0)^2}\\
\stackrel{\rm a.s.}{\longrightarrow}
\frac{\bm_n^\top\bSigma^{-1}_n\bm_n \bm_0^\top\bSigma^{-1}_n\bm_0-  (\bm_n^\top\bSigma^{-1}_n\bm_0)^2}
{(c+ \bm_n^\top\bSigma^{-1}_n\bm_n) \bm_0^\top\bSigma^{-1}_n\bm_0-(\bm_n^\top\bSigma^{-1}_n\bm_0)^2}=\alpha^*
$$
and, similarly,
\begin{equation*}
\hat{\beta}^*=(1-\hat{\alpha}^*)\, \dfrac{\bby_n^\top\bS^{-1}_n\bm_0}{\bm_0^\top\bS^{-1}_n\bm_0}\stackrel{\rm a.s.}{\longrightarrow} (1-\alpha^*)\, \dfrac{\bm_n^\top{\bSigma}^{-1}_n\bm_0}{\bm_0^\top{\bSigma}^{-1}_n\bm_0}=\beta^*
\end{equation*}
$\text{for}~ p/n\to  c \in (0, \infty)~~\text{as} ~n\rightarrow\infty$. Thus the theorem is proved. \hfill $\Box$

\bigskip
The proof of Theorem~\ref{th4} is based on Lemma~\ref{lem3}. Let
$\shatnull =\bbx_n^\top \bS^{-1}_n\bbx_n- {(\bm_0^\top \bS^{-1}_n\bbx_n)^2}/{\bm_0^\top \bS^{-1}_n\bm_0}$ and $
\Rhatnull = {\bm_0^\top \bS^{-1}_n\bbx_n}/{\bm_0^\top \bS^{-1}_n\bm_0}$.

\begin{lemma}\label{lem3}
Let $\mathbf{y}_1,\ldots,\mathbf{y}_n$ be a random sample of independent vectors such that $\mathbf{y}_i\sim\mathcal{N}_p(\bm_n,\bSigma_n)$ for all $i \in \{1,\ldots,n\}$. Then for any $p$ and $n$ with $n>p$, the following statements hold.
\begin{enumerate}[(a)]
\item $ {n(n-p+1)} \shatnull/\{{(n-1)(p-1)}\} \sim F_{p-1,n-p+1,n\,\snull}$,  i.e., the non-central $F$-distribution with $p-1$ and $n-p+1$ degrees of freedom and non-centrality parameter $\snull$.
\item $\Rhatnull | \shatnull=y \sim \mathcal{N}\left[ \Rnull , (1+\frac{n}{n-1}y)/\{n \bm_0^\top \bSigma^{-1}_n\bm_0\}\right]$.
\end{enumerate}
\end{lemma}

\vspace{0.3cm}
\noindent{\it Proof  of Lemma~\ref{lem3}.}
The results are obtained following the proof of Theorem 3.1 in \cite{bodnarschmid2008}, where similar statements are presented in the case $\bm_0=\bi_p$. \hfill $\Box$

\bigskip
\noindent{\it Proof  of Theorem~\ref{th4}.}
Let $\shatnullc =(1-c)\shatnull -c$.
Then
$
\sqrt{n}\left(\shatnullc -\snull\right)=\sqrt{n}\, \{(1-c)\shatnull - c-\snull\}
=\sqrt{n}\, \{(1-c)\shatnull/c-1- \snull/c\})c.
$
An application of Lemma~\ref{lem3}.a and Lemma 3 in \cite{bodnar_reiss_2016} leads to
\begin{equation}\label{sig_s2}
\sqrt{n} (\shatnullc-\snull )\rightsquigarrow  \mathcal{N} (0,\sigma_s^2 ) \quad \text{with} \quad \sigma_s^2= 2\left(c+2\snull\right)+\frac{2}{1-c}\left(c+\snull\right)^2.
\end{equation}
From Theorem~\ref{th3}, we get
$\hat{\alpha}^*=1-c(1-c)^{-1}/\shatnull =1-{c}/{(c+\shatnullc)}$
and
\begin{align*}
\hat{\beta}^*=\frac{c}{c+\shatnullc }\Rhatnull \stackrel{d}{=}\frac{c}{c+\shatnullc }\left(\Rnull +\frac{\sqrt{1+\frac{n}{n-1}\frac{\shatnullc +c}{1-c}}}
{\sqrt{\bm_0^\top \bSigma^{-1}_n\bm_0}}z_0\right),
\end{align*}
where
\begin{align*}
\sqrt{n}\left(\begin{array}{c}
  \shatnullc-\snull \\
  z_0\\
          \end{array}
        \right)\rightsquigarrow  \mathcal{N}\left[\left(
                                                   \begin{array}{c}
                                                     0 \\
                                                     0 \\
                                                   \end{array}
                                                 \right),
                                                 \left(
                                                   \begin{array}{cc}
                                                     \sigma_s^2  & 0 \\
                                                     0 & 1\\
                                                   \end{array}
                                                 \right)
        \right].
\end{align*}
with $\sigma_s^2$ as in \eqref{sig_s2}. Now, an application of the Delta method (see, e.g., Theorem 3.7 in \cite{dasgupta2008}) leads to
\begin{align*}
\sqrt{n}\left(\begin{array}{c}
  \hat{\alpha}^*-\alpha^* \\
  \hat{\beta}^*-\beta^*\\
          \end{array}
        \right)\rightsquigarrow  \mathcal{N}\left[ \left(
                                                   \begin{array}{c}
                                                     0 \\
                                                     0 \\
                                                   \end{array}
                                                 \right),
                                                 \bOmega
        \right]
\end{align*}
with
\begin{multline*}
\bOmega=
\left(\begin{array}{cc} \frac{c}{(c+\snull)^2} & 0 \\ \frac{c}{(c+\snull)^2}\Rnull  & \frac{c}{c+\snull}\frac{\sqrt{1+\frac{\snull+c}{1-c}}}{\sqrt{\bm_0^\top \bSigma^{-1}_n\bm_0}} \\\end{array}\right)
\left(\begin{array}{cc}\sigma_s^2  & 0 \\0 & 1\\\end{array}\right)
\left(\begin{array}{cc} \frac{c}{(c+\snull)^2} &\frac{c}{(c+\snull)^2}\Rnull  \\ 0 & \frac{c}{c+\snull}\frac{\sqrt{1+\frac{\snull+c}{1-c}}}{\sqrt{\bm_0^\top \bSigma^{-1}_n\bm_0}} \\\end{array}\right) \\
= \left(\begin{array}{cc} \frac{c^2\sigma_s^2}{(c+\snull)^4} &\frac{c^2 \sigma_s^2}{(c+\snull)^4}\Rnull  \\
\frac{c^2 \sigma_s^2}{(c+\snull)^4}\Rnull  &\frac{c^2 \sigma_s^2}{(c+\snull)^4}\Rnull ^2 + \frac{c^2}{(c+\snull)^2}\frac{1+\frac{\snull+c}{1-c}}{\bm_0^\top \bSigma^{-1}_n\bm_0} \\\end{array}\right).
\end{multline*}
Hence the theorem is proved. \hfill $\Box$

\bigskip
In the proof of Theorem~\ref{th4} we make use of the following lemma.

\begin{lemma}\label{lem4}
Assume (A1)--(A2) and let the elements of $\bx_n$ possess uniformly bounded $4+\varepsilon$, $\varepsilon>0$, moments. Let $\boldsymbol{\theta}$ and $\boldsymbol{\xi}$ be the universal nonrandom vectors from the set $\{\bm_0, \bm_n \}$. Then
\begin{align*}
   p^{-\gamma}\left|\boldsymbol{\xi}^\top\bS_n^{-}\boldsymbol{\theta}- c^{-1}(c-1)^{-1} \boldsymbol{\xi}^\top \bSigma^{-1}_n\boldsymbol{\theta}\right| \stackrel{\rm a.s.}{\longrightarrow} 0 ,\quad
   \bbx_n^\top\bSigma_n^{1/2}\bS_n^{-}\bSigma_n^{1/2}\bbx_n \stackrel{\rm a.s.}{\longrightarrow} \frac{1}{c-1}  ,\quad
 p^{-\gamma/2}\bbx_n^\top\bSigma_n^{1/2}\bS_n^{-}\boldsymbol{\theta} \stackrel{\rm a.s.}{\to }0
\end{align*}
for $p/n\to  c \in (1, \infty)$ as $n\rightarrow\infty$, where $\bbx_n=n^{-1}\bx_n \bi_n$ stands for the sample mean vector calculated from $\bx_n$.
\end{lemma}

\noindent{\it Proof  of Lemma~\ref{lem4}.}
For $\boldsymbol{\theta} \in  \{\bm_0, \bm_n \}$, one has
$\left\|\bSigma^{-1/2}_n\boldsymbol{\theta}\right\|^2\le \lambda_{\max}(\bSigma^{-1}_n)\left\|\boldsymbol{\theta}\right\|^2<\infty$
following assumptions (A1) and (A2). The rest of the proof follows from Lemma 5.6 in \cite{bodnar2016}. \hfill $\Box$

\bigskip
\noindent{\it Proof  of Theorem~\ref{th5}.}
For $\boldsymbol{\theta},\bxi \in  \{\bm_0, \bm_n \}$, the application of Lemma~\ref{lem4} leads to
$$
p^{-\gamma} |\bm_0^\top\bS_n^{-}\bm_0- c^{-1}(c-1)^{-1} \bm_0^\top \bSigma^{-1}_n\bm_0 | \stackrel{\rm a.s.}{\longrightarrow} 0,
$$
$$
p^{-\gamma} |\bm_0^\top\bS_n^{-}\bby_n-c^{-1}(c-1)^{-1} \bm_0^\top \bSigma^{-1}_n\bm_n |
\le p^{-\gamma} |\bm_0^\top\bS_n^{-}\bSigma_n^{1/2}\bbx_n |
+ p^{-\gamma} |\bm_0^\top\bS_n^{-}\bm_n-c^{-1}(c-1)^{-1}\bm_0^\top \bSigma^{-1}_n\bm_n |  \stackrel{\rm a.s.}{\longrightarrow} 0,
$$
and
\begin{multline*}
 p^{-\gamma} |\bby_n^\top\bS_n^{-}\bby_n-c^{-1}(c-1)^{-1} \bm_n^\top \bSigma^{-1}_n\bm_n-(c-1)^{-1} |
\le p^{-\gamma} |\bm_n^\top\bS_n^{-}\bm_n-c^{-1}(c-1)^{-1} \bm_n^\top \bSigma^{-1}_n\bm_n |\\
+p^{-\gamma} |\bbx_n^\top\bSigma_n^{1/2}\bS_n^{-}\bSigma_n^{1/2}\bbx_n-(c-1)^{-1} |
+2p^{-\gamma} |\bm_n^\top\bS_n^{-}\bSigma_n^{1/2}\bbx_n |\stackrel{\rm a.s.}{\longrightarrow} 0.
\end{multline*}

Hence,
\begin{align*}
\hat{\alpha}^*= \frac{ \{\bby_n^\top \bS^{-}_n\bby_n -(c-1)^{-1}  \}\bm_0^\top \bS^{-}_n\bm_0-(\bby_n^\top \bS^{-}_n\bm_0)^2}{\bby_n^\top \bS^{-}_n\bby_n \bm_0^\top \bS^{-}_n\bm_0-(\bby_n^\top \bS^{-}_n\bm_0)^2}
\stackrel{\rm a.s.}{\longrightarrow}
\frac{\bm_n^\top\bSigma^{-1}_n\bm_n \bm_0^\top\bSigma^{-1}_n\bm_0-  (\bm_n^\top\bSigma^{-1}_n\bm_0)^2}
{(c+ \bm_n^\top\bSigma^{-1}_n\bm_n) \bm_0^\top\bSigma^{-1}_n\bm_0-(\bm_n^\top\bSigma^{-1}_n\bm_0)^2}=\alpha^*\end{align*}
and, similarly,
\begin{equation*}
\hat{\beta}^*=(1-\hat{\alpha}^*) \, \dfrac{\bby_n^\top\bS^{-}_n\bm_0}{\bm_0^\top\bS^{-}_n\bm_0}\stackrel{\rm a.s.}{\longrightarrow} (1-\alpha^*) \, \dfrac{\bm_n^\top{\bSigma}^{-1}_n\bm_0}{\bm_0^\top{\bSigma}^{-1}_n\bm_0}=\beta^*.
\end{equation*}
$\text{for}~ p/n\to  c \in (1, \infty)~~\text{as} ~n\rightarrow\infty$. This completes the proof of the theorem. \hfill $\Box$

\section*{Acknowledgments}

We are grateful to Professors Christian Genest and Philippe Vieu as well as two anonymous reviewers for their helpful suggestions. We also gratefully acknowledge the comments from the participants at the Fourth International Workshop on Functional and Operatorial Statistics (A Coru\~{n}a), the European Meeting of Statisticians (Amsterdam), the WISE-CASE Workshop on Econometrics and Statistics (Xiamen, China), and the Joint Meeting of the German Mathematical Society and the Polish Mathematical Society (Poznan). Finally, we thank Cheng Wang for providing us the MATLAB code for his estimation method.

\end{document}